\documentclass[11pt, leqno]{article}

\usepackage[margin=25mm]{geometry}
\usepackage{amsmath, amssymb, amsthm, amscd, mathrsfs, mathtools, extarrows}
\usepackage[all]{xy}
\usepackage{enumerate}
\usepackage{url}

\makeatletter
\@addtoreset{equation}{section}
\makeatother

\theoremstyle{theorem}
\newtheorem{theorem}{Theorem}[subsection]
\newtheorem{definition}[theorem]{Definition}
\newtheorem{proposition}[theorem]{Proposition}
\newtheorem{lemma}[theorem]{Lemma}
\newtheorem{corollary}[theorem]{Corollary}

\newtheorem{hypothesis}[theorem]{Hypothesis}

\theoremstyle{definition}
\newtheorem{remark}[theorem]{Remark}{\rm}
\newtheorem{example}[theorem]{Example}{\rm}
{\rm}
{\rm}
\newtheorem{_empty}[theorem]{}{\rm}
{\rm}

\newcommand{\Cs}{\mathscr{C}}
\newcommand{\Cone}{{\rm Cone}}
\newcommand{\Coh}{{\rm Coh}}
\newcommand{\Bf}{\mathfrak{B}}
\newcommand{\D}{\mathcal{D}}
\newcommand{\Db}{\mathbb{D}}

\newcommand{\DVR}{\mathrm{DVR}}

\newcommand{\E}{\mathcal{E}}
\newcommand{\Ed}{\mathcal{E}^\dag}

\newcommand{\F}{\mathcal{F}}

\newcommand{\G}{\mathcal{G}}
\newcommand{\Homs}{\mathcal{H}\kern -.5pt om}
\newcommand{\id}{{\rm id}}
\newcommand{\I}{\mathcal{I}}

\newcommand{\Isoc}{{\rm Isoc}^\dag}

\newcommand{\Lb}{\mathbb{L}}
\newcommand{\Lotimes}{\otimes^{\Lb}}
\newcommand{\LD}[1]{\underrightarrow{LD}^{b}_{\Qb, \mathrm{#1}}}

\newcommand{\MIC}{{\rm MIC}}

\renewcommand{\O}{\mathcal{O}}
\newcommand{\oh}{\mathrm{oh}}
\newcommand{\ovhol}{\mathrm{ovhol}}
\newcommand{\oc}{\mathrm{oc}}
\newcommand{\ovcoh}{\mathrm{ovcoh}}
\newcommand{\op}{\mathrm{op}}
\newcommand{\Pf}{\mathfrak{P}}

\newcommand{\Qb}{\mathbb{Q}}
\newcommand{\Qf}{\mathfrak{Q}}
\newcommand{\R}{\mathcal{R}}
\newcommand{\Rb}{\mathbb{R}}
\newcommand{\RG}{\mathbb{R}\underline{\Gamma}}
\newcommand{\rig}{\mathrm{rig}}
\newcommand{\red}{{\rm red}}
\newcommand{\simto}{\xrightarrow{\sim}}

\newcommand{\Strat}{{\rm Strat}}

\newcommand{\Tf}{\mathfrak{T}}

\newcommand{\Us}{\mathscr{U}}
\newcommand{\Uf}{\mathfrak{U}}

\newcommand{\V}{\mathcal{V}}

\newcommand{\Vs}{\mathscr{V}}
\newcommand{\W}{\mathcal{W}}

\newcommand{\Xf}{\mathfrak{X}}

\newcommand{\Yf}{\mathfrak{Y}}

\newcommand{\Zf}{\mathfrak{Z}}
\renewcommand{\sp}{{\rm sp}}
\newcommand{\Spa}{{\rm Spa}\,}

\newcommand{\Spec}{{\rm Spec}\,}
\newcommand{\Spf}{{\rm Spf}\,}

\newcommand{\ul}{\underline}
\newcommand{\ol}{\overline}
\newcommand{\wh}{\widehat}
\newcommand{\wt}{\widetilde}

\title{Overholonomicity of Overconvergent F-isocrystals over Laurent Series Fields}
\author{Yuanmin Liu}
\date{}

\begin{document}
	\maketitle
	\begin{abstract}
		We prove that $K$-valued overconvergent $F$-isocrystals over smooth $k((t))$-varieties are overholonomic and construct the triangulated category of overholonomic complexes of Frobenius type which is stable under six operations and contains overconvergent $F$-isocrsytals.
	\end{abstract}
	\tableofcontents
	
	\section{Introduction}
	Let $\V$ be a complete discrete valuation ring of mixed characteristic $(0, p)$ with perfect residue field $k$, let $K$ be its field of fractions and let $\sigma : \V \to \V$ be a lift of the $s$-th Frobenius map on $k$. 
	
	For a scheme $X$ over $k$ of finite type, Berthelot defined the $p$-adic Weil cohomology called rigid cohomology $H^*_\rig(X/K)$. It is a variation of de Rham cohomology over $p$-adic rigid analytic spaces with suitable condition of convergence. The coefficients are called overconvergent $F$-isocrystals, which are coherent modules with integrable connections equipped with Frobenius structures. However, the category of overconvergent $F$-isocrystals is insufficient; they are $p$-adic entities of lisse $l$-adic sheaves (or Weil sheaves), which are not stable under Grothendieck's six operations. Berthelot tried to understand these operations via arithmetic $\D$-modules (\cite{berthelotDI} \cite{berthelotDII}), in which he defined six operations and the Frobenius pullback functor. There should be a $p$-adic Riemann-Hilbert correspondence giving an equivalence between ``holonomic arithmetic $\D$-modules'' and the category of ``constructible isocrystals'' stable under six operations. This is achieved by the effort of many mathematicians and now is fully available \cite[Theorem 8.0.1]{lazda2023rigidification}.
	
	Let us move on to the case of $k((t))$-varieties. In this setting, the Amice ring $\E_K$ is a complete discrete valuation field of mixed characteristic $(0, p)$ with residue field $k((t))$ and one gets the rigid cohomology $H^*_{\rig}(X/\E_K)$. These cohomology groups have the additional structure of $(\sigma, \nabla)$-modules, which have a relation with $p$-adic representations. From this viewpoint, coming from geometry, these $(\sigma, \nabla)$-modules must have some sort of finite monodromy and it is natural to expect that they are defined over the Robba ring $\R_K$; the $p$-adic monodromy theorem proved by Andr\'{e}-Mebkhout-Kedlaya \cite{andre2002HasseArf}\cite{mebkhout2002analogue}\cite{kedlaya2004p} states that $(\sigma, \nabla)$-modules over $\R_K$ are quasi-unipotent. The problem is that two rings $\R_K$ and $\E_K$ have no inclusions. Lazda-P\'{a}l resolved this problem by refining Berthelot's definition to get rigid cohomology valued over the bounded Robba ring $\Ed_K$, which is the intersection of $\E_K$ and $\R_K$.
	This is done by taking the compactification not into proper $k((t))$-schemes but into proper $k[[t]]$-schemes. On the arithmetic $\D$-module side, Caro defined them over formal schemes locally having finite $p$-bases over $\V$ (Definition \ref{def:p-basis}). This allows one to consider arithmetic $\D$-modules over strictly semistable formal schemes over $\V[[t]]$ (cf. Definition \ref{def:ssp}), which is quite useful for local computation as we will see.
	
	Having these two objects, overconvergent isocrystals and arithmetic $\D$-modules, it is natural to expect that results in the classical case are still valid; in particular there should be a category of arithmetic $\D$-modules stable under six operations and containing overconvergent $F$-isocrystals. In order to prove this, one needs overholonomicity of overconvergent $F$-isocrystals proved in \cite{caro2012overholonomicity}. Overholonomic $\D$-modules is abstractly defined to form the largest triangulated full subcategory stable under proper pushforward, pullbacks, duals and direct summands. For overholonomicity, one needs a structure theorem of overconvergent $F$-isocrystals, namely the semistable reduction theorem proved by Kedlaya in the series of papers \cite{kedlaya2007semistable} \cite{kedlaya2008semistable} \cite{kedlaya2009semistable} \cite{kedlaya2011semistable}. The author proved the semistable reduction theorem for overconvergent $F$-isocrystals in the sense of Lazda-P\'al in  \cite{liu2026semistable} and we apply it to the study of arithmetic $\D$-modules in this paper.
	
	\subsection*{Main Results}
	We construct the category $D^b_{\oh, F}(X/K)$ of overholonomic arithmetic $\D$-modules of Frobenius type for any realizable scheme $X$ over $k[[t]]$ (see \ref{emp:realizable}). Examples of realizable schemes over $k[[t]]$ are quasi-projective $k((t))$-schemes and quasi-projective $k$-schemes. In Definition \ref{def:six}, we also define six operations
	$$ f^!, \, f^+ : D^b_{\oh, F}(X/K) \to D^b_{\oh, F}(X'/K) $$
	$$ f_!, \, f_+ : D^b_{\oh, F}(X'/K) \to D^b_{\oh, F}(X/K) $$
	for a morphism $f : X' \to X$ of realizable schemes over $k[[t]]$ and
	$$ - \otimes_X - : D^b_{\oh, F}(X/K) \times D^b_{\oh, F}(X/K) \to D^b_{\oh, F}(X/K)$$
	$$ \Db_X : D^b_{\oh, F}(X/K)^\op \to D^b_{\oh, F}(X/K). $$
	For a realizable scheme $X$, we have the category $\Isoc_F(X/K)$ of overconvergent isocrystals of Frobenius type. It is equivalent to the one defined by Berthelot when $X$ is quasi-projective over $k$ and equivalent to the one defined by Lazda-P\'al when $X$ is quasi-projective over $k((t))$. We also define the fully faithful functor $\sp_{X, +} : \Isoc(X/K) \to D^b_{\oh, F}(X/K)$ with following compatibilities of operations:
	$$ \sp_{X, +}(E^\vee) \simeq \Db_X \sp_{X, +}(E),$$
	$$ \sp_{X, +}(E \otimes F) \simeq \sp_{X, +}(E) \, \wt{\otimes} \, \sp_{X, +}(F)[\delta_X]$$
	and
	$$ \sp_{X', +}f^*_K(E) \simeq f^! \sp_{X, +}(E)[\delta_{X'} - \delta_X]$$ 
	for $E, F \in \Isoc_F(X/K)$ and a morphism $f : X' \to X$ of realizable $k[[t]]$-schemes (Theorem \ref{thm:sp-compatibility}).
	
	\subsection*{Structure of the paper}
	In section 2, we define overholonomic $\D$-modules and prove stability properties by operations except tensor products.
	
	In section 3, we define the functor $\sp_+$ from the category of overconvergent isocrystals to the category of arithmetic $\D$-modules. We also define the category of arithmetic $\D$-modules ``devissable in isocrystals''. This is necessary for proving the stability of overholonomic $F$-complexes under tensor products.
	
	In section 4, we prove the comparison theorem between log rigid cohomology and rigid cohomology, following \cite{caro2012overholonomicity}. We are only able to prove it under stronger condition about exponents, but it will be sufficient for our purpose.
	
	In section 5, we translate the comparison theorem of section 4 in terms of log arithmetic $\D$-modules. As an important consequence, we prove that for a strictly semistable pair $(X, Z)$ over $k[[t]]$, an overconvergent isocrystal over $(X-Z, X)$ coming from a log convergent isocrystal with nilpotent residues is overholonomic.
	
	In section 6, we define the category $D^b_{\oh, F}$ and six operations, based on the results from previous sections.
	
	\subsection*{Acknowledgement}
	The author is thankful to his advisor Atsushi Shiho for helpful discussions. He is also thankful to Daniel Caro, who gave useful comments on this paper and supported the author during the visit to Caen university.
	
	This work was supported by the World-leading INnovative Graduate Study Program for Frontiers of Mathematical Sciences and Physics, The University of Tokyo.
	
	\subsection*{Notations and conventions}
	\begin{itemize}
		\item $\V$ is a complete discrete valuation ring of mixed characteristic $(0, p)$ with perfect residue field $k$. Let $K$ be its field of fractions and let $\sigma : \V \to \V$ be a lift of the $s$-th Frobenius map on $k$. 
		\item A formal scheme is $p$-adic and locally noetherian. We say that a formal scheme $\Pf$ is of (locally) formal finite type over $\V$ if it is (locally) of finite type over some $\Db^d_\V := \Spf \V[[t]] \times_\V \cdots \times_\V \Spf \V[[t]]$ (for varying $d$).
		\item A formal scheme is denoted by fraktur alphabets such as $\Pf$, and we often denote by the corresponding usual alphabets such as $P$ to denote its special fiber unless otherwise stated.
		\item A log structure on a scheme or on a formal scheme is always assumed to be fine.
	\end{itemize}
	
	\section{Overholonomicity of $\D$-modules}
	In this section, we define the notion of formal schemes over $\V$ locally having finite $p$-basis, over which Caro defined arithmetic $\D$-modules and their overholonomicity, the condition which we are going to prove that overconvergent $F$-isocrystals should satisfy. The notations and results of this section are mostly cited from \cite{caro2021arithmetic}.
	
	\subsection{Finite $p$-bases}
	\begin{definition}[{\cite[1.1.1]{caro2021arithmetic}}]
		A morphism of $k$-schemes $f : X \to Y$ is said to be relatively perfect if the following diagram is cartesian:
		\[\xymatrix{
			X \ar[d]_{f} \ar[r]^{{\rm Fr}_X} & X \ar[d]^{f} \\
			Y \ar[r]^{{\rm Fr}_Y} & Y
		}\]
	\end{definition}
	\begin{proposition}[{\cite[1.2]{kato1991explicit}}]
		Let $f : X \to Y$ be a morphism of $k$-schemes. Consider the following conditions:
		\begin{enumerate}[(i)]
			\item $f$ is etale;
			\item $f$ is relatively perfect;
			\item $f$ is formally etale.
		\end{enumerate}
		We have the implication (i) $\Rightarrow$ (ii) $\Rightarrow$ (iii).
	\end{proposition}
	
	\begin{definition}[cf. {\cite[Definition 1.1.7]{caro2021arithmetic}}]\label{def:p-basis}
		\begin{enumerate}[(a)]
			\item Let $f : X \to Y$ be a morphism of $k$-schemes. $t_1, \dots, t_n \in \Gamma(X, \O_X)$ is said to be a finite $p$-basis of $f$ or $X$ over $Y$ when $f$ is clear from the context if the induced morphism $X \to \mathbb{A}^n_Y$ is relatively perfect.
			\item Let $f : \Xf \to \Yf$ be a morphism of formal $\V$-schemes. We say that $f$ is relatively perfect if $f$ is formally etale and $f_0 : X \to Y$ is relatively perfect.
			\item Let $f : \Xf \to \Yf$ be a morphism of formal $\V$-schemes. $t_1, \dots, t_n \in \Gamma(\Xf, \O_\Xf)$ is said to be a finite $p$-basis of $f$ or $\Xf$ over $\Yf$ when $f$ is clear from the context if the induced morphism $\Xf \to \wh{\mathbb{A}}^n_\Yf$ is relatively perfect.
		\end{enumerate}
		The cardinality of finite $p$-bases coincides with the rank of $\Omega^1_{X/Y}$, which we will denote by $\delta_{X/Y}.$ It is a function assigning to each connected components of $X$ a nonnegative integer.
	\end{definition}
	
	\begin{example}
		\begin{enumerate}[(i)]
			\item $t_1, \cdots, t_n$ is a finite $p$-basis of $\Spec k[[t]][t_1, \dots, t_n]/(t - t_1 \cdots t_r)$ over $\Spec k$ and is that of $\Spf \V[[t]]\langle t_1, \dots, t_n \rangle/(t - t_1 \cdots t_r)$ over $\Spf \V.$
			\item A smooth morphism $f : X \to Y$ of $k$-schemes locally has finite $p$-bases. An etale coordinate is finite $p$-basis.
			\item A smooth morphism $f : \Xf \to \Yf$ of formal $\V$-schemes locally has finite $p$-bases. An etale coordinate is finite $p$-basis.
		\end{enumerate}
	\end{example}
	
	\subsection{Overholonomic $\D$-modules}
	\begin{definition}
		Let $\DVR(\V)$ be the category of local homomorphisms $\V \to \W$ to complete discrete valuation rings $\W$ of mixed characteristic $(0, p)$ with perfect residue fields. 
		A data of coefficient $\Cs$ is a datum assigning each object $\W$ of $\DVR(\V)$ and a connected formal $\W[[t]]$-scheme $\Pf$ of finite type, locally having finite $p$-bases over $\W$ to a strict full subcategory $\Cs(\Pf/\W)$ of $\LD{coh}(\wh{\D}_{\Pf/\W}^{(\bullet)})$.
	\end{definition}
	
	\begin{remark}
		The definition is none other than the {\it restricted data of absolute coefficients over $\V$} defined in \cite[Definition 13.1.4]{caro2021arithmetic}. Also note that $\LD{coh}(\wh{\D}_{\Pf/\W}^{(\bullet)})$ is canonically equivalent to $D^b_{\mathrm{coh}}(\D^\dag_{\Pf/\W, \Qb})$ via taking colimits, \cite[Proposition 5.2.8]{caro2021arithmetic}.
	\end{remark}
	
	\begin{_empty}
		For an object $\W$ of $\DVR(\V)$ and a formal $\W[[t]]$-scheme $\Pf$ of finite type, locally having finite $p$-bases, we have the following operations:
		\begin{itemize}
			\item The local cohomological functor $\RG^\dag_{Z} : \LD{\mathrm{qc}}(\wh{\D}_{\Pf/\W}^{(\bullet)}) \to \LD{\mathrm{qc}}(\wh{\D}_{\Pf/\W}^{(\bullet)})$ for a closed subscheme $Z \subseteq \Pf$.
			\item The localization outside a divisor $(^\dag D) : \LD{\mathrm{qc}}(\wh{\D}_{\Pf/\W}^{(\bullet)}) \to \LD{\mathrm{qc}}(\wh{\D}_{\Pf/\W}^{(\bullet)}).$
			\item The dual functor $\Db : \LD{\mathrm{coh}}(\wh{\D}_{\Pf/\W}^{(\bullet)})^\op \to \LD{\mathrm{coh}}(\wh{\D}_{\Pf/\W}^{(\bullet)}).$
		\end{itemize}
		For a morphism $f : \Pf \to \Qf$ of formal $\W[[t]]$-schemes of finite type, locally having finite $p$-bases, we have the following operations:
		\begin{itemize}
			\item The extraordinary pullback $f^! : \LD{\mathrm{qc}}(\wh{\D}_{\Qf/\W}^{(\bullet)}) \to \LD{\mathrm{qc}}(\wh{\D}_{\Pf/\W}^{(\bullet)})$
			\item The pushforward $f_+ : \LD{\mathrm{qc}}(\wh{\D}_{\Pf/\W}^{(\bullet)}) \to \LD{\mathrm{qc}}(\wh{\D}_{\Qf/\W}^{(\bullet)}).$
		\end{itemize}
		Note that we omit the superscript $(\bullet)$ on the functors for simplifying the notations.
	\end{_empty}
	
	\begin{definition}
		Let $\W$ be an object of $\DVR(\V)$, $\Pf$ a connected formal $\W[[t]]$-scheme of finite type, locally having finite $p$-bases over $\W$ and $n \geq 0$, $r \geq -1$ two integers. We define that a complex $\E \in \LD{\mathrm{qc}}(\wh{\D}_{\Pf/\W}^{(\bullet)})$ is $(n, r)$-overholonomic by induction on $r$ as follows:
		\begin{itemize}
			\item $\E$ is $(n, -1)$-overholonomic if $\delta_P < n$ or if $\delta_P \leq n$ and $\E \in \LD{\mathrm{coh}}(\wh{\D}_{\Pf/\W}^{(\bullet)})$.
			\item $\E$ is $(n, 0)$-overholonomic if $\delta_P < n$ or if $\E$ is $(n, -1)$-overholonomic and for any smooth morphism $f : \Pf' \to \Pf$ from a connected $\W[[t]]$-formal scheme $\Pf'$ with $\delta_{P'} \leq n$ and any divisor $T'$ of $P'$, we have $(^\dag T') f^! \E \in \LD{\mathrm{coh}}(\wh{\D}_{\Pf'/\W}^{(\bullet)})$.
			\item $\E$ is $(n, r)$-overholonomic for $r > 0$ if $\delta_P < n$ or if $\E$ is $(n, r-1)$-overholonomic and for any smooth morphism $f : \Pf' \to \Pf$ from a connected $\W[[t]]$-formal scheme $\Pf'$ with $\delta_{P'} \leq n$ and any divisor $T'$ of $P'$, $\Db (^\dag T') f^! \E$ is $(n, r-1)$-overholonomic.
		\end{itemize}
		We say that $\E$ is $r$-overholonomic if $\E$ is $(n, r)$-overholonomic for all $n \geq 0$ and that $\E$ is overholonomic if $\E$ is $r$-overholonomic for all $r \geq -1$. When $\Pf$ is not necessarily connected, we say that $\E$ is $r$-overholonomic (resp. overholonomic) if it is over each connected component of $\Pf$. The full subcategory of overholonomic complexes over $\Pf$ is denoted by $\LD{ovhol}(\wh{\D}^{(\bullet)}_{\Pf})$ or simply by $\LD{ovhol}(\Pf)$. The essential image of $\LD{ovhol}(\wh{\D}^{(\bullet)}_{\Pf})$ via $\varinjlim : \LD{coh}(\wh{\D}^{(\bullet)}_{\Pf}) \simto D^b_{\rm coh}(\D^\dag_{\Pf, \Qb})$ is denoted by $ D^b_{\rm ovhol}(\D^\dag_{\Pf, \Qb})$.
	\end{definition}
	
	\begin{remark}
		$(n, r)$-overholonomicity is the analogue of \cite[Definition 2.3.5]{caro2012overholonomicity} with a slight modification to clarify the definition in the case of $\delta_P < n,$ which is implicitly used in the proof of \cite[Theorem 2.3.13]{caro2012overholonomicity}. The definition here is a priori more restrictive than the that of \cite[Definition 13.1.2]{caro2012overholonomicity}, where the extraordinary pullbacks are only tested for the projections $f : \Pf \times \Qf \to \Pf$ with $\Qf$ a smooth formal $\V$-schemes, but our definition is equivalent nevertheless.
		
		Also note that this definition contains the classical case where $\Pf$ is a smooth formal $\V$-schemes, over which the definition is equivalent to the classical one.
	\end{remark}
	
	\begin{proposition} \label{prop:nr-ovhol}
		Let $\W$ be an object of $\DVR(\V)$, $n \geq 0$ and $r \geq -1$ two integers.
		\begin{enumerate}[(i)]
			\item Let $g : \Qf \to \Pf$ be a projective morphism of connected formal $\W[[t]]$-schemes of finite type, locally having finite $p$-bases over $\W$. For $d = \delta_Q - \delta_{P}$, if $\F$ is an $(n+d, r)$-overholonomic complex over $\Qf$, then $g_+\F$ is $(n, r)$-overholonomic.
			\item Let $g : \Qf \to \Pf$ be a smooth morphism of connected formal $\W[[t]]$-schemes of finite type, locally having finite $p$-bases over $\W$. If $\E$ is an $(n, r)$-overholonomic complex over $\Pf$, then $g^!\E$ is $(n, r)$-overholonomic.
		\end{enumerate}
	\end{proposition}
	\begin{proof}
		(i): We may assume $\delta_X \leq n$. When $r = -1$, this is simply the fact that pushforward by projective morphisms preserves coherence, \cite[Proposition 10.7.2]{caro2012overholonomicity}. Suppose $r \geq 0$. Let $f : \Pf' \to \Pf$ be a smooth morphism from a connected $\W[[t]]$-formal scheme $\Pf'$ with $\delta_{P'} \leq n$ and let $T'$ be a divisor of $P'$. Consider the following cartesian square:
		\[\xymatrix{
			\Qf' \ar[r]^{g'} \ar[d]_{f'} & \Pf' \ar[d]^{f} \\
			\Qf \ar[r]^{g} & \Pf.
		}\]
		Using the base change theorem \cite[Corollary 10.6.11]{caro2021arithmetic} and the commutation of localizations with pushforwards \cite[Corollary 7.4.3]{caro2021arithmetic}, we get $(^\dag T')f^!g_+\F \simeq (^\dag T')g'_+f'^! \F \simeq g'_+ (^\dag g'^{-1}(T')) f'^! \F$. Since $\F$ is $(n+d, 0)$-overholonomic, $(^\dag g'^{-1}(T')) f'^! \F$ both belongs to $\LD{\mathrm{coh}}(\wh{\D}_{\Qf'/\W}^{(\bullet)}).$ Again using the fact that $g'_+$ preserves coherence, $(^\dag T')f^!g_+\F$ belongs to $\LD{\mathrm{coh}}(\wh{\D}_{\Pf'/\W}^{(\bullet)})$ and we conclude that $g_+\F$ is $(n, 0)$-overholonomic. Lastly when $r > 0$, the proof goes as the case where $r = 0$; the only difference is to use the relative duality $g'_+ \Db \simeq \Db g'_+$ proved in \cite[Theorem 10.7.3]{caro2021arithmetic}.
		
		(ii): When $r = -1$, the assertion follows from \cite[Proposition 7.5.2]{caro2021arithmetic}. We leave the proof to readers.
	\end{proof}
	
	\begin{_empty}
		Let $f : \W' \to \W$ be a morphism of $\DVR(\V)$, let $\Pf$ be a formal $\W[[t]]$-scheme of finite type, locally having finite $p$-bases over $\W$ and let $f : \Pf' := \Pf \times_{\Spf \W} \Spf \W' \to \Pf$ be the projection morphism by abuse of notation. Then we have the base change functor
		$$f^* : \LD{qc}(\wh{\D}^{(\bullet)}_{\Pf/\W}) \to \LD{qc}(\wh{\D}^{(\bullet)}_{\Pf'/\W'}), \ \E^{(\bullet)} \mapsto \E^{(\bullet)} \wh{\otimes}^{\mathbb{L}}_{\W} \W'.$$
		We can also define the functor
		$$f^* : D(\D^\dag_{\Pf/\W, \Qb}) \to D(\D^\dag_{\Pf'/\W', \Qb}), \ \E \mapsto \D^\dag_{\Pf'/\W', \Qb} \otimes^{\mathbb{L}}_{f^{-1}\D^\dag_{\Pf/\W, \Qb}} f^{-1}\E.$$
		Note that these two functors are exact, they map coherent complexes to coherent complexes and they are compatible under the equivalences $\varinjlim : \LD{qc}(\wh{\D}^{(\bullet)}_{\Pf/\W}) \simto D^b_{\rm coh}(\D^\dag_{\Pf/\W, \Qb})$ and that with primes. Also one easily checks that base change functors are compatible with extraordinary pullbacks, pushforwards, dual functors and local cohomological functors.
	\end{_empty}
	
	\begin{definition}
		Let $\W$ be an object of $\DVR(\V)$, $\Pf$ a connected formal $\W[[t]]$-scheme of finite type, locally having finite $p$-bases over $\W$ and $n \geq 0$, $r \geq -1$ two integers. We say that a complex $\E \in \LD{\mathrm{qc}}(\wh{\D}_{\Pf/\W}^{(\bullet)})$ is overholonomic after any base change if for any morphism $f : \W \to \W'$ of $\DVR(\V)$ of the form $\W \to \W \otimes_{W(l)} W(l')$ for $l$ the residue field of $\W$ and $l'$ a perfect field extending $l$, the base change $f^* \E$ is overholonomic. We denote by $\LD{\oh}(\wh{\D}_{\Pf/\W}^{(\bullet)})$ the full subcategory of overholonomic complexes after any base change and by $D^b_{\oh}(\D^\dag_{\Pf/W, \Qb})$ its essential image under $\varinjlim : \LD{\mathrm{coh}}(\wh{\D}_{\Pf/\W}^{(\bullet)}) \to D^b_{\mathrm{coh}}(\D^\dag_{\Pf/W, \Qb})$.
	\end{definition}
	
	\begin{proposition} \label{prop:ovhol}
		Let $\W$ be an object of $\DVR(\V)$, $\Pf$ a formal $\W[[t]]$-scheme locally with finite $p$-basis over $\W$. Let $?$ be either $\ovhol$ or $\oh$.
		\begin{enumerate}[(i)]
			\item $\O_{\Pf}^{(\bullet)}$ is overholonomic after any base change;
			\item If $\E \in \LD{?}(\wh{\D}^{(\bullet)}_{\Pf})$ and $\E' \in \LD{qc}(\wh{\D}^{(\bullet)}_{\Pf})$ its direct summand, then $\E' \in \LD{?}(\wh{\D}^{(\bullet)}_{\Pf})$;
			\item If $\E \in \LD{?}(\wh{\D}^{(\bullet)}_{\Pf})$, then $\Db(\E) \in \LD{?}(\wh{\D}^{(\bullet)}_{\Pf})$;
			\item For an exact triangle $\E' \to \E \to \E'' \to$ in $\LD{qc}(\wh{\D}^{(\bullet)}_{\Pf})$, if two of them are overholonomic (after any base change), then so is the other;
			\item For $\E \in \LD{qc}(\wh{\D}^{(\bullet)}_{\Pf})$, if there is an open cover of $\Pf$ over each of which $\E$ is overholonomic (after any base change), then $\E$ is overholonomic (after any base change);
			\item If $\E \in \LD{?}(\wh{\D}^{(\bullet)}_{\Pf})$, then for any divisor $T$ of $P$ and any closed subscheme $Z$ of $P$, the complexes $(^\dag T)\E$ and $\RG^\dag_{Z} \E$ are in $\LD{?}(\wh{\D}^{(\bullet)}_{\Pf})$;
			\item If $\E \in \LD{?}(\wh{\D}^{(\bullet)}_{\Pf})$ and $f : \Pf \to \Qf$ is a projective morphism to a formal $\W[[t]]$-scheme $\Qf$ locally with finite $p$-basis over $\W$, then $f_+ \E \in \LD{?}(\wh{\D}^{(\bullet)}_{\Qf})$;
			\item If $\E \in \LD{?}(\wh{\D}^{(\bullet)}_{\Pf})$ and $g : \Pf' \to \Pf$ is a quasi-projective morphism from a formal $\W[[t]]$-scheme $\Pf'$ locally with finite $p$-basis over $\W$, then $g^! \E \in \LD{?}(\wh{\D}^{(\bullet)}_{\Pf'})$.
		\end{enumerate}
	\end{proposition}
	\begin{proof}
		These are the consequences of \cite[Theorem 13.4.7]{caro2021arithmetic} and \cite[Example 13.4.8]{caro2021arithmetic}.
	\end{proof}
	
	\section{Devissage in isocrystals}
	\subsection{Isocrystals}
	In this subsection, we give a slightly generalized definition of $K$-valued overconvergent isocrystals in order to allow realization over formal schemes locally of formal finite type.
	\begin{definition}
		A formal frame is a triple $(X, Y, \Pf)$ of a formal $\V$-scheme $\Pf$ locally of formal finite type, a closed subscheme $Y$ of $\Pf$ and an open subscheme $X$ of $Y$.
		A morphism of formal frames is a triple $(f, g, u) : (X', Y', \Pf') \to (X, Y, \Pf)$ of morphisms such that
		\[\xymatrix{
			X' \ar[r] \ar[d]_{f} & Y' \ar[r] \ar[d]_{g} & \Pf' \ar[d]_{u} \\
			X \ar[r] & Y \ar[r] & \Pf \\
		}
		\]
		commutes.
		
		We say that a formal frame is formally smooth if $\Pf$ is locally smooth over some  $\Db^d_\V := \Spf \V[[t]] \times_\V \cdots \times_\V \Spf \V[[t]]$ in an neighborhood of $X$ (for varying $d$).
	\end{definition}
	
	\begin{definition}
		Let $(X, Y, \Pf)$ be a formal frame. Then we have the specialization map $\sp : \Pf_K \to \Pf$ and the map $[-] : \Pf_K \to \Pf_K$ sending the point to its maximal generalization. We denote by $[\sp]$ the composition of these two maps. For a subset $S$ of $P$, we denote by $]S[_{\Pf}$ the inverse image of $S$ by $[\sp]$. Let $i : ]X[_\Pf \hookrightarrow ]Y[_\Pf$ be the inclusion and put $\O_{]X[_\Pf} := i^{-1}\O_{]Y[_\Pf}$ and $j^\dag_X (-) = i_* i^{-1}(-).$  Note that $[\sp]$ is anti-continuous; hence $]Y[_\Pf$ is an adic space locally of formal finite type over $K$ and $]X[_\Pf$ has the structure of a germ.
		
		For a morphism of frames $(f, g, u) : (X', Y', \Pf') \to (X, Y, \Pf),$ it induces a morphism of germs  $]X'[_{\Pf'} \to ]X[_{\Pf}$, which we will denote by $]u[$ {\rm (cf. \cite[\S 2]{abe2025proper})}.
	\end{definition}
	
	\begin{definition}
		Let $(X, Y, \Pf)$ be a formally smooth formal frame.
		\begin{itemize}
			\item The category of overconvergent isocrystals $\Isoc(X, Y, \Pf/K)$ is defined in the following way: an object is a collection of coherent $\O_{]X[_\Pf}$-modules (or equivalently a coherent $j^\dag_{X}\O_{]Y[_{\Pf}}$-module) $\E_{\Pf'}$ for formal frames $(X', Y', \Pf') $ over $(X, Y, \Pf)$ and a collection of isomorphisms $\epsilon_{u} : ]u[^*(\E_{\Pf'}) \to \E_{\Pf''}$ for morphisms of formal frames $(f, g, u) : (X'', Y'', \Pf'') \to (X', Y', \Pf')$ satisfying the cocycle condition. A morphism of overconvergent isocrystals is a collection of morphisms $\E_{\Pf} \to \E'_{\Pf}$ commuting with isomorphisms $\epsilon_u$ and $\epsilon'_u$.
			\item The category of overconvergent stratified modules $\Strat^\dag(X, Y, \Pf/K)$ is defined in the following way: an object is a datum $(\E, \epsilon)$ of a coherent $\O_{]X[_\Pf}$-module $\E$ and an isomorphism $\epsilon : p_1^*\E \to p_0^*\E$ satisfying the cocycle condition $]p_{02}[^*\epsilon = (]p_{01}[^*\epsilon) \circ (]p_{12}[^*\epsilon)$, where $p_i : (X, Y, \Pf \times_\V \Pf) \to (X, Y, \Pf)$ and $p_{ij} : (X, Y, \Pf \times_\V \Pf \times_\V \Pf) \to (X, Y, \Pf \times_\V \Pf)$ are natural projection maps. A morphism of overconvergent stratifications is a morphism $\E \to \E'$ commuting with isomorphisms $\epsilon$ and $\epsilon'$.
			\item The category of stratified modules $\Strat^\dag(X, Y, \Pf/K)$ is defined in the following way: an object is a datum $(\E, \epsilon)$ of a coherent $\O_{]X[_\Pf}$-module $\E$ and a collection of compatible isomorphisms $\epsilon^{(n)} : ]p_1^{(n)}[^*\E \to ]p_0^{(n)}[^*\E$ satisfying $\epsilon^{(0)} = \id$ and the cocycle condition $p^*_{02}\epsilon^{(n)} = (p^*_{01}\epsilon^{(n)}) \circ (p^*_{12}\epsilon^{(n)})$, where $p^{(n)}_i : (X, Y, \Pf^{(n)}) \to (X, Y, \Pf)$ are projection maps from the $n$-th infinitesimal neighborhood $\Pf^{(n)}$ of $\Pf$ in $\Pf \times_\V \Pf$ and $p_{ij} : (X, Y, \Pf_2^{(n)}) \to (X, Y, \Pf^{(n)})$ are natural projection maps from  the $n$-th infinitesimal neighborhood $\Pf_2^{(n)}$ of $\Pf$ in $\Pf \times_\V \Pf \times_\V \Pf$. A morphism of overconvergent stratifications is a morphism $\E \to \E'$ commuting with isomorphisms $\epsilon$ and $\epsilon'$.
			\item The category of coherent modules with integrable connections $\MIC(X, Y, \Pf/K)$ is defined in the following way: an object is a datum $(\E, \nabla)$ of a coherent $\O_{]X[_\Pf}$-module $\E$ and a connection $\E \to \Omega^1_{]X[_\Pf/K}$ which is integrable. A morphism of coherent modules with integrable connections is a morphism $\E \to \E'$ commuting with $\nabla$ and $\nabla'$.
			\item Given a morphism of formally smooth formal frames $\theta = (f, g, u) : (X', Y', \Pf') \to (X, Y, \Pf)$, we have the pullback functors $\theta^*$ between these categories. We sometimes simply denote it by $u^*_K$ if $\theta$ is clear from the context.
		\end{itemize}
	\end{definition}
	
	\begin{proposition}
		Let $(X, Y, \Pf)$ be a formally smooth formal frame. Then there is a sequence of functors
		$$\Isoc(X, Y, \Pf/K) \to \Strat^\dag(X, Y, \Pf/K) \to \Strat(X, Y, \Pf/K) \to \MIC(X, Y, \Pf),$$
		whose left-most and right-most arrows are equivalences and the middle one is fully faithful.
	\end{proposition}
	\begin{proof}
		The proof is quite similar to those of classical case. Let us prove that the middle arrow is fully faithful. The proof is almost the same as in \cite[Proposition 5.13]{lazda2016rigid}. Following the argument, it reduces to prove that for $S = \V[[t_1]] \wh{\otimes}_\V \cdots \wh{\otimes}_\V \V[[t_d]] = \Gamma(\Db^d_{\V}, \O)$ and the ring $\O_m$ of sections on the closed tube $[\Db^d_k]^m_{\Db^d_\V \times_\V \Db^d_\V}$ of radius $p^{-1/m}$, the ring map
		$$S_K \wh{\otimes}_K S_K \to S_K[[\tau_1, \dots, \tau_d]], \ f \otimes g \mapsto f(t_1, \dots, t_d) g(t_1 + \tau_1, \dots, t_d + \tau_d)$$
		induces the injective map $\O_m \to S_K[[\tau_1, \dots, \tau_d]]$. Since we have
		$$\O_m = (S_K \wh{\otimes}_K S_K) \langle T_1, \dots, T_d \rangle / (pT_1 - (1 \otimes t_1 - t_1 \otimes 1)^m, \dots, pT_d - (1 \otimes t_d - t_d \otimes 1)^m),$$
		we get the assertion by an easy calculation.
	\end{proof}
	
	\begin{definition}
		Let $(X, Y, \Pf)$ be a formally smooth formal frame. The essential image of $\Strat^\dag(X, Y, \Pf/K)$ in $\MIC(X, Y, \Pf/K)$ is denoted by $\MIC^\dag(X, Y, \Pf/K)$ and is called the category of overconvergent modules with integrable connections.
	\end{definition}
	
	\begin{proposition}
		Let $(X, Y, \Pf)$ be a formally smooth formal frame and $\E \in \MIC(X, Y, \Pf/K)$. Then $\E$ is overconvergent if and only if $\Pf$ can be covered by affine open formal subschemes smooth over some  $\Db^d_\V$ in an neighborhood of $X$ such that for a basis $dx_1, \dots, dx_n$ of $\Omega^1_{\Pf/\V}$ around $X$ and any $\eta \in (0, 1)$, $\E$ comes from some coherent module $\E_V$ with integrable connection over an affinoid open neighborhood $V$ of $]X[_{\Pf}$ in $\Pf_K$ and $\left\| \displaystyle \frac{\ul{\partial}^{\ul{k}}}{\ul{k} !} e \right\| \eta^{\ul{k}} \to 0$ as $\ul{k} \to \infty$ for all $e \in \Gamma(V, \E_V)$, where $\| \cdot \|$ is a Banach norm on $\Gamma(V, \E_V)$.
	\end{proposition}
	\begin{proof}
		The proof is similar to the classical case \cite[Lemme 2.2.12]{berthelot1996cohomologie}. See also \cite[Proposition 2.68]{lazda2016rigid}.
	\end{proof}
	
	\begin{proposition}[Weak fibration theorem]
		Let $(id, g, u) : (X, Y', \Pf') \to (X, Y, \Pf)$ be a morphism of formal frames with $g$ proper and $u$ etale in an open neighborhood of $X$. Then we have an isomorphism of germs $]u[ \,:\, ]X[_{\Pf'} \to ]X[_{\Pf}$.
	\end{proposition}
	\begin{proof}
		Similarly as in \cite[Proposition 1.3.1]{berthelot1996cohomologie}, there are isomorphisms of an open neighborhood of $]X[_{\Pf'}$ and that of $]X[_{\Pf}$, which preserves the closed subspaces $]X[_{\Pf'}$ and $]X[_{\Pf}$.
	\end{proof}
	
	\begin{proposition}
		Let $(id, g, u) : (X, Y', \Pf') \to (X, Y, \Pf)$ be a morphism of formally smooth formal frames with $g$ proper and $u$ formally smooth in an open neighborhood of $X$. Then the pullback functor $u^* : \Isoc(X, Y, \Pf / K) \to \Isoc(X, Y', \Pf' / K)$ is an equivalence of categories.
	\end{proposition}
	\begin{proof}
		Similarly as in \cite[Th\'eor\`em 2.3.1, Th\'eor\`em 2.3.5]{berthelot1996cohomologie}.
	\end{proof}
	
	\begin{remark}
		By these propositions, when $(X, Y, \Pf)$ is a formal frame such that $X$ and $Y$ are of finite type over $k$ and $\Pf$ is of finite type over $\V$, the category $\Isoc(X, Y, \Pf/K)$ defined here coincides with the usual one. The same also holds for $X$ of finite type over $k((t))$, $Y$ of finite type over $k[[t]]$ and $\Pf$ of finite type over $\V[[t]]$.
	\end{remark}
	
	\subsection{Partially overconvergent $\D$-modules}
	\begin{_empty}[Overcoherent modules with overconvergent singularity along $Z$]
		Let $\Pf$ be a formal scheme over $\V[[t]]$ of finite type, locally having finite $p$-bases over $\V$. For a closed subscheme $Z$ of $P$, define the functor $(^\dag Z) : \LD{qc}(\wh{\D}_{\Pf}^{(\bullet)}) \to \LD{qc}(\wh{\D}_{\Pf}^{(\bullet)})$ by the exact triangle
		$$ \RG_{Z}^\dag \E \to \E \to (^\dag Z) \to.$$
		When $Z$ is a divisor, it is simply the localization functor outside a divisor. We denote by $\LD{qc}(\wh{\D}_{\Pf}^{(\bullet)}(Z))$ the full subcategory of $\LD{qc}(\wh{\D}_{\Pf}^{(\bullet)})$ consisting of objects $\E$ such that $\E \to (^\dag Z) \E$ is an isomorphism, or equivalently $\RG^\dag_Z \E = 0$. We denote by $\LD{ovcoh}(\wh{\D}_{\Pf}^{(\bullet)}(Z))$ the full subcategory of $\LD{qc}(\wh{\D}_{\Pf}^{(\bullet)}(Z))$ consisting of objects $\E$ such that for any smooth morphism $f : \Pf' \to \Pf$ and any divisors $T \subseteq T'$ of $P'$ containing $f^{-1}(Z)$, $(^\dag T)f^!\E$ belongs to $(^\dag T') \E \in \LD{coh}(\wh{\D}_{\Pf'}^{(\bullet)}(T))$. Also we let $\LD{oc}(\wh{\D}_{\Pf}^{(\bullet)}(Z))$ the full subcategory of $\LD{qc}(\wh{\D}_{\Pf}^{(\bullet)}(Z))$ consisting of objects $\E$ such that for any morphism $f : \V \to \W = \V \otimes_{W(k)} W(l)$ in $\DVR(\V)$ for some perfect field $l$ over $k$, $f^*\E$ belongs to $\LD{ovcoh}(\wh{\D}_{\Pf'}^{(\bullet)}(Z'))$ where $\Pf' = \Pf \times_{\Spf \V} \Spf \W$ and $Z' = f^{-1}(Z).$
	\end{_empty}
	
	\begin{lemma}
		Let $?$ be either $\ovcoh$ or $\oc.$ We have the following:
		\begin{itemize}
			\item $\LD{?}(\wh{\D}_{\Pf}^{(\bullet)}(Z))$ is a triangulated full subcategory of $\LD{qc}(\wh{\D}_{\Pf}^{(\bullet)}(Z))$.
			\item For $\E \in \LD{qc}(\wh{\D}_{\Pf}^{(\bullet)}(Z))$, $\E$ belongs to $\LD{?}(\wh{\D}_{\Pf}^{(\bullet)}(Z))$ is a Zariski local condition on $\Pf$. 
			\item For a smooth morphism $f : \Pf' \to \Pf$ and $\E \in \LD{?}(\wh{\D}_{\Pf}^{(\bullet)}(Z))$, $f^! \E$ belongs to  $\LD{qc}(\wh{\D}_{\Pf'}^{(\bullet)}(f^{-1}(Z)))$.
			\item For a closed subscheme $Z'$ of $P$, the functor $(^\dag Z')$ maps $\LD{?}(\wh{\D}_{\Pf}^{(\bullet)}(Z))$ into $\LD{?}(\wh{\D}_{\Pf}^{(\bullet)}(Z)) \cap \LD{?}(\wh{\D}_{\Pf}^{(\bullet)}(Z \cup Z'))$.
			\item For a locally closed subscheme $Y'$ of $P$, the functor $\RG^\dag_{Y'}$ maps $\LD{?}(\wh{\D}_{\Pf}^{(\bullet)}(Z))$ into $\LD{?}(\wh{\D}_{\Pf}^{(\bullet)}(Z))$.
		\end{itemize}
	\end{lemma}
	\begin{proof}
		Similar to \cite[15.3.7]{caro2025arith}.
	\end{proof}
	
	\begin{proposition}
		Let $\E \in \LD{ovcoh}(\wh{\D}_{\Pf}^{(\bullet)}(Z))$ and let $\Uf = \Xf - Z$. If $\E|_\Uf = 0$, then $\E = 0$. 
	\end{proposition}
	\begin{proof}
		When $Z$ is a divisor, then $\E$ is a coherent $\wh{\D}_{\Pf}^{(\bullet)}(Z)$-module. Then by the faithful flatness \cite[Theorem 3.3.4]{caro2021arithmetic}, one gets the assertion. The general case goes by induction on the number of divisors $T_1, \dots, T_r$ such that $Z = T_1 \cap \cdots \cap T_r$, see \cite[Proposition 15.3.7.6]{caro2025arith}.
	\end{proof}
	
	\begin{proposition}
		Let $?$ be either $\ovcoh$ or $\oc.$ For $f : \Pf' \to \Pf$ a proper morphism of formal $\V[[t]]$-schemes locally having finite $p$-bases over $\V$ and $Z' = f^{-1}(Z)$, the functor $f_+$ maps $\LD{?}(\wh{\D}_{\Pf'}^{(\bullet)}(f^{-1}(Z)))$ into $\LD{?}(\wh{\D}_{\Pf}^{(\bullet)}(Z))$.
	\end{proposition}
	\begin{proof}
		Similar to \cite[Proposition 15.3.7.8]{caro2025arith}.
	\end{proof}
	
	\begin{_empty}
		Let $\Pf$ be a formal scheme over $\V[[t]]$ of finite type, $Z$ and $X$ two closed subschemes of $P$. Let $\LD{povcoh}(X, \Pf, Z / \V)$ denote the full subcategory of $\LD{ovcoh}(\wh{\D}_{\Pf}^{(\bullet)}(Z))$ consisting of objects $\E$ such that $\E|_{\Pf - X} = 0$ (or equivalently $\RG^\dag_{X} \E = 0$, see \cite[Proposition 15.3.8.2]{caro2025arith}). Define $\LD{poc}(X, \Pf, Z/\V)$ similarly. Let $?$ be either ${\rm povcoh}$ or ${\rm poc}$; we define the canonical $t$-structure on $\LD{?}(X, \Pf, Z/\V)$ as follows: $\E \in \LD{?}(X, \Pf, Z/\V)$ belongs to $\underrightarrow{LD}^{b, \lesseqqgtr 0}_{\Qb, ?}(X, \Pf, Z / \V)$ if and only if $\E|_{\Pf - Z} \in \underrightarrow{LD}^{b, \lesseqqgtr 0}_{\Qb, \mathrm{qc}}(\wh{\D}^{(\bullet)}_{\Pf-Z/\V})$. We denote by $\tau_Z$ and $H^n_Z$ the associated truncation functors and cohomology functors. For the proof that this indeed defines a $t$-structure, see \cite[Proposition 15.3.8.12]{caro2025arith}.
		
		Remark that in general the $t$-structure on $\LD{?}(X, \Pf, Z / \V)$ is not in general induced from that of $\LD{qc}(\wh{\D}_{\Pf}^{(\bullet)})$ unless $Z$ is a divisor, see \cite[Example 15.3.8.13]{caro2025arith}.
	\end{_empty}
	
	\begin{proposition}
		Let $?$ be either $\ovcoh$ or $\oc.$ For $f : \Pf' \to \Pf$ a morphism of formal $\V[[t]]$-schemes locally having finite $p$-bases over $\V$ and $Z' = f^{-1}(Z)$, the functor $f^!$ maps $\LD{?}(\wh{\D}_{\Pf}^{(\bullet)}(Z))$ into $\LD{?}(\wh{\D}_{\Pf'}^{(\bullet)}(f^{-1}(Z)))$.
	\end{proposition}
	\begin{proof}
		Similar to \cite[Proposition 15.3.8.27]{caro2025arith}.
	\end{proof}
	
	\subsection{$\sp_+$ functor, divisorial case}
	\begin{_empty} \label{emp:MIC-ddag}
		Let $(Y, X, \Pf)$ be a formal frame with $\Pf$ of finite type over $\V[[t]]$, locally having finite $p$-bases over $\V$ and smooth over $\V[[t]]$ in an open neighborhood of $X$. Suppose that there is a divisor $T$ of $P$ such that $Y = X - T$.
		
		We define the category $\MIC^\ddag(X, \Pf, T / K)$ to be the category of coherent $\D^\dag_\Pf(^\dag T)_\Qb$-modules $\E$ supported on $X$, which means that the canonical morphism $\RG^\dag_{X} \E \simto \E$ is an isomorphism, such that for any open subset $\Uf$ of $\Pf$ and a formal $\V[[t]]$-scheme $\Xf$ locally having finite $p$-bases over $\V$ lifting $X \cap U$ and any closed immersion $\iota : \Xf \hookrightarrow \Uf$ (over $\V$) lifting $X \cap U \hookrightarrow U$, the coherent $\D^\dag_{\Xf}(^\dag T \cap X \cap U)_\Qb$-module $\iota_\alpha^!(\E|_{\Uf})[-\delta_{X \cap U} + \delta_U]$ is coherent over $\O_\Uf(^\dag T \cap U)_\Qb.$
		
		Note that for a coherent $\D^\dag_\Pf(^\dag T)_\Qb$-module $\E$ supported on $X$, to check that $\E$ belongs to $\MIC^\ddag(X, \Pf, T / K)$, one only needs to check coherence by some lift $\iota_\alpha : \Xf_\alpha \hookrightarrow \Uf_\alpha$ for an open cover $\{ \Uf_\alpha \}$ of $\Pf.$ This is checked as follows : given two lifts $\iota : \Xf_ \hookrightarrow \Uf_\alpha$ and $\iota' : \Xf \hookrightarrow \Uf_\alpha$ there is a canonical isomorphism (called the glueing isomorphism) $\tau_{\iota, \iota'} : \iota^! \simto \iota'$ by \cite[Proposition 8.5.1]{caro2021arithmetic} and it is independent of the choice of $\iota$; for another lift $\Xf'$ of $X \cap \Uf_\alpha$, since it is isomorphic to $\Xf$ over $\V$ by formal smoothness, it is also independent of the choice of lifts.
		
		Let $f : (Y', P', \Pf') \to (Y, P, \Pf)$ be a morphism of formal frames and suppose that $\Pf'$ is of finite type over $\V[[t]]$, locally having finite $p$-bases over $\V$ and smooth over $\V[[t]]$ in an open neighborhood of $X'$. Suppose that there is a divisor $T'$ of $P'$ such that $Y' = X' - T'$ and $f^{-1}(T) \subseteq T'$. Then the functor $f^* := f_{T, T'}^![-\delta_{P'}+\delta_{P}]$ induces a functor $\MIC^\ddag(P \Pf, T / K) \to \MIC^\ddag(P', \Pf', T' / K)$ as in \cite[Corollary 11.2.3.5]{caro2025arith} and $\Db_T$ induces a functor $\MIC^\ddag(P, \Pf, T / K)^\op \to \MIC^\ddag(X, \Pf, T / K)$ as in \cite[11.2.6.3]{caro2025arith}; the arguments work as well in our case since the local descriptions are the same.
	\end{_empty}
	
	\begin{lemma} \label{lem:sp_*}
		In the situation of {\rm \ref{emp:MIC-ddag}}, we have quasi-inverses $$\sp_* : \MIC^\dag(Y, P, \Pf / K) \leftrightarrows \MIC^\ddag(P, \Pf, T / K) : \sp^*.$$ Furthermore, we have natural isomorphisms
		$$f^* \sp_* E \simeq \sp_* f_K^* E,$$
		$$ \Db_T (\sp_* E) \simeq \sp_*(E^\vee) $$
		and
		$$ \sp_*(E) \otimes^{\dag, \mathbb{L}}_{\O_{\Pf, \Qb}} \sp_*(F) \simeq \sp_*(E \otimes F) $$
		for any $E, F \in \MIC^\dag(Y, P, \Pf / K)$. Moreover, the natural isomorphism $f^*\sp_* \simeq \sp_* f^*_K$ is compatible with glueing isomorphisms $\epsilon_{f, f'} : f'^*_K \simto f^*_K$ and $\tau_{f, f'} : f'^* \simto f^*$.
	\end{lemma}
	\begin{proof}
		We can argue similarly as in the classical case \cite[11.2.7.1, Proposition 11.2.7.3, Corollary 11.2.7.5]{caro2025arith}.
	\end{proof}
	
	\begin{_empty}
		Let $\{\Pf_\alpha\}$ be an open covering of $\Pf$, $\Pf_{\alpha_1 \cdots \alpha_i} := \Pf_{\alpha_1} \cap \cdots \cap \Pf_{\alpha_i}$, $X_{\alpha_1 \dots \alpha_i} := X \cap P_{\alpha_1 \cdots \alpha_i}$, $T_{\alpha_1 \cdots \alpha_i} := T \cap P_{\alpha_1 \cdots \alpha_i}$. Let $\Xf_{\alpha_1 \cdots \alpha_i}$ be an affine formal $\V[[t]]$-scheme locally with finite $p$-bases over $\V$ lifting $X_{\alpha_1 \cdots \alpha_i}$, let $\iota_{\alpha_1 \cdots \alpha_i} : \Xf_{\alpha_1 \cdots \alpha_i} \hookrightarrow \Pf_{\alpha_1 \cdots \alpha_i}$ be a closed immersion lifting $X_{\alpha_1 \cdots \alpha_i} \hookrightarrow P_{\alpha_1 \cdots \alpha_i}$ and let $p_{j_1 \cdots j_r}^{\alpha_1 \cdots \alpha_i} : \Xf_{\alpha_1 \cdots \alpha_i} \to \Xf_{\alpha_{j_1} \cdots \alpha_{j_r}}$ be a lift of the open immersion $X_{\alpha_1 \cdots \alpha_i} \hookrightarrow X_{\alpha_{j_1} \cdots \alpha_{j_r}}$. Note that $\iota$ or $p$ are not necessarily $\V[[t]]$-morphisms; they are $\V$-morphisms. We call the data $(\Xf_{\alpha_1}, \Xf_{\alpha_1 \alpha_2}, \Xf_{\alpha_1 \alpha_2 \alpha_3}, \iota_{\alpha_1}, \iota_{\alpha_1 \alpha_2}, \iota_{\alpha_1 \alpha_2 \alpha_3}, p_{i}^{\alpha_1 \alpha_2}, p_{i}^{\alpha_1 \alpha_2 \alpha_3}, p_{ij}^{\alpha_1 \alpha_2 \alpha_3})$ a lifting system of $X \hookrightarrow \Pf$ with respect to the open cover $\{\Pf_\alpha\}$. For simplicity we denote the lifting system by $(\Xf_\alpha).$
		
		Given a lifting system $(\Xf_\alpha)$, we the category $\Coh((\Xf_\alpha), T / K)$ as follows: an object is a data $(\E_\alpha, \theta_{\alpha_1 \alpha_2})$ of a collection of coherent $\D^\dag_{\Xf_\alpha}(^\dag T_\alpha)_\Qb$-modules $\E_\alpha$ and a collection of isomorphisms $\theta^{\alpha_1 \alpha_2} : p_{2}^{\alpha_1 \alpha_2, !} \E_{\alpha_2} \simto p_{1}^{\alpha_1 \alpha_2, !} \E_{\alpha_1}$ satisfying the cocycle condition: let $\theta^{\alpha_1 \alpha_2 \alpha_3}_{i j} : p_{j}^{\alpha_1 \alpha_2 \alpha_3, !} \E_j \simto p_{i}^{\alpha_1 \alpha_2 \alpha_3, !} \E_i$ be the isomorphism making the diagram
		\[\xymatrix{
			p^{\alpha_1 \alpha_2 \alpha_3, !}_{ij} p^{\alpha_i \alpha_j, !}_{2}\E_j \ar[r]^-{\tau}_-{\sim} \ar[d]^-{\sim}_-{p^{\alpha_1 \alpha_2 \alpha_3, !}_{ij} \theta^{\alpha_i \alpha_j}} & p_{j}^{\alpha_1 \alpha_2 \alpha_3, !} \E_j \ar[d]^{\theta^{\alpha_1 \alpha_2 \alpha_3}_{i j}}\\
			p^{\alpha_1 \alpha_2 \alpha_3, !}_{ij} p^{\alpha_i \alpha_j, !}_{1}\E_i \ar[r]^-{\tau}_-{\sim} & p_{i}^{\alpha_1 \alpha_2 \alpha_3, !} \E_i;
		}\]
		then we have $\theta^{\alpha_1 \alpha_2 \alpha_3}_{1 2} \theta^{\alpha_1 \alpha_2 \alpha_3}_{2 3} = \theta^{\alpha_1 \alpha_2 \alpha_3}_{1 3}$. A morphism of $(\E_\alpha, \theta_{\alpha_1 \alpha_2}) \to (\E'_\alpha, \theta'_{\alpha_1 \alpha_2})$ is defined in a natural way.
		
		Let $\MIC^\ddag(X, \Pf, T / K)$ be the full subcategory of $\Coh((\Xf_\alpha), T / K)$ whose objects are those $(\E_\alpha, \theta_{\alpha_1 \alpha_2})$ with $\E_\alpha$ coherent over $\O_{\Xf}(^\dag X_\alpha \cap T_\alpha)_\Qb$.
		Let $\Coh(X, \Pf, T / K)$ be the category of coherent $\D^\dag_{\Pf}(^\dag T)_\Qb$-modules supported over $X$.
	\end{_empty}
	
	\begin{lemma}
		There is a natural equivalence $u_0^! : \Coh(X, \Pf, T / K) \to \Coh((\Xf_\alpha), T / K), \ \E \mapsto (\iota_\alpha^! (\E|_{\Pf_\alpha}), \star)$. The essential image of $\MIC^\ddag(X, \Pf, T / K)$ is $\MIC^\ddag((\Xf_\alpha), T / K)$.
	\end{lemma}
	\begin{proof}
		Similarly as in \cite[Lemma 5.3.6]{caro2019arithmetic}.
	\end{proof}
	
	\begin{_empty}
		Let $(\Xf_\alpha)$ be a lifting system. We can define the category $\MIC^\dag((\Xf_\alpha), T/K)$ similarly as $\MIC^\ddag((\Xf_\alpha), T/K)$ and define the equivalence $u_K^* : \MIC^\dag(X, \Pf, T / K) \simto \MIC^\dag((\Xf_\alpha), T/K)$. Note that the equivalence of $u_K^*$ comes from Lemma \ref{lem:indep-frame}. We have the canonical functor $\sp_* : \MIC^\dag((\Xf_\alpha), T / K) \to \MIC^\ddag((\Xf_\alpha), T / K)$, $(E_\alpha, \theta^{\alpha_1 \alpha_2}) \mapsto (\sp_*E_\alpha, \sp_*\theta^{\alpha_1 \alpha_2}).$
		
		We define the functor $\sp_{X \hookrightarrow \Pf, T, +} : \MIC^\dag(X, \Pf, T / K) \to \MIC^\ddag(X, \Pf, T / K)$ by making the following diagram commutative:
		\[\xymatrix@C=3.9pc{
			\MIC^\dag(X, \Pf, T / K) \ar[r]^-{\sp_{X \hookrightarrow \Pf, T, +}} \ar[d]^{u_K^*} & \MIC^\ddag(X, \Pf, T / K) \ar[d]^{u_0^!}\\
			\MIC^\dag((\Xf_\alpha), T / K) \ar[r]^-{\sp_*} & \MIC^\ddag((\Xf_\alpha), T / K).\\
		}\]
		Note that $\sp_{X \hookrightarrow \Pf, T, +}$ is an equivalence of categories since the other three in the diagram are equivalences.
		It is straightforward to check that the functor $\sp_{X \hookrightarrow \Pf, T, +}$ is independent of the choice of lifting systems $(\Xf_\alpha)$.
	\end{_empty}
	
	\begin{proposition} \label{prop:sp+-pullback}
		In the situation of {\rm \ref{emp:MIC-ddag}}, we have a commutative diagram
		\[\xymatrix@C=3.9pc{
			\MIC^\dag(X', \Pf', T' / K) \ar[r]^{\sp_{X' \hookrightarrow \Pf', T', +}} \ar[d]_{f_K^*} & \MIC^\ddag(X, \Pf', T' / K) \ar[d]_{f^*} \\
			\MIC^\dag(X, \Pf, T / K) \ar[r]^{\sp_{X \hookrightarrow \Pf, T, +}} & \MIC^\ddag(X, \Pf, T / K) \\
		}\]
	\end{proposition}
	\begin{proof}
		Choose a lifting system of $(\Xf_\alpha)$ of $X \subseteq \Pf$ and a lifting system $(\Xf'_\beta)$ of $X' \subseteq \Pf'$. After refinement, we may assume that they are indexed by a same set and there exist lifts $f_{\alpha_1 \cdots \alpha_n} : \Xf'_{\alpha_1 \cdots \alpha_n} \to \Xf_{\alpha_1 \cdots \alpha_n}$. Then using Lemma \ref{lem:sp_*} for each $f_{\alpha_1 \cdots \alpha_n}$, from the definition of $\sp_+$, we get the assertion. 
	\end{proof}
	
	\begin{proposition}
		In the situation of {\rm \ref{emp:MIC-ddag}}, for any $E \in \MIC^\dag(X, \Pf, T/K)$, we have canonical isomorphism
		$$ \sp_{X \hookrightarrow \Pf, T, +}(E^\vee) \simeq \Db_T(\sp_{X \hookrightarrow \Pf, T, +} E).$$
	\end{proposition}
	\begin{proof}
		The assertion is local and use Lemma \ref{lem:sp_*} and the fact that the pushforward functor $i_+$ by a closed immersion commutes with the dual functors $\Db$ \cite[Corollary 8.4.4]{caro2021arithmetic}.
	\end{proof}
	
	\begin{proposition}
		In the situation of {\rm \ref{emp:MIC-ddag}}, for any $E, F \in \MIC^\dag(X, \Pf, T/K)$, we have canonical isomorphism
		$$ \sp_{X \hookrightarrow \Pf, T, +}(E \otimes F) \simeq \sp_{X \hookrightarrow \Pf, T, +}(E) \otimes^{\dag, \mathbb{L}}_{\O_{\Pf, \Qb}} \sp_{X \hookrightarrow \Pf, T, +}(F)[\delta_X - \delta_P].$$
	\end{proposition}
	\begin{proof}
		The assertion is local and we suppose that we have a lift $i : \Xf \hookrightarrow \Pf$. Let $\E = \sp_{X \hookrightarrow \Pf, T, +}(E)$ and $\F = \sp_{X \hookrightarrow \Pf, T, +}(F)$. By \cite[7.2.3]{caro2021arithmetic}, we have $i^! \E \otimes^{\dag, \mathbb{L}}_{\O_{\Xf, \Qb}} i^! \F \simeq i^! (\E \otimes^{\dag, \mathbb{L}}_{\O_{\Pf, \Qb}} \F) [\delta_X - \delta_P].$ Since $i^! \E = \sp_* E$, $i^! \F = \sp_* F$ by definition and $i^! \E \otimes^{\dag, \mathbb{L}}_{\O_{\Xf, \Qb}} i^! \F = \sp_* (E \otimes F)$ by Lemma \ref{lem:sp_*}, applying $i_+$ we get the assertion.
	\end{proof}
	
	\begin{proposition} \label{prop:dir-sum-isoc}
		Let
		\[\xymatrix{
			X' \ar@{^(->}[r] \ar[d]^{a} & \Pf' \ar[d]^{f}\\
			X \ar@{^(->}[r] & \Pf
		}\]
		be a commutative diagram where $X$ and $X'$ are $k[[t]]$-schemes locally having finite $p$-bases over $k$, $\Pf$ and $\Pf'$ are formal $\V[[t]]$-schemes locally having finite $p$-bases over $\V$, horizontal arrows are closed immersions, and $a$ a flat morphism and $f$ a smooth projective morphism. Let $T$ be a divisor of $P$, $T' := f^{-1}(T)$ and suppose that $T \cap X$ is a divisor of $X$ (hence $T' \cap Y$ is a divisor of $Y$). Then for a coherent $\D^\dag_{\Pf}(^\dag T)_\Qb$-module $\E$ in the essential image of $\sp_{X \hookrightarrow \Pf, T, +}$, $\E$ is a direct summand of $f_+\RG_{X'}^\dag f^!\E$.
	\end{proposition}
	\begin{proof}
		Similarly as in \cite[Proposition 9.6.5]{caro2021arithmetic}, using the fact that $\Db_{T'}$ locally corresponds to the dual of isocrystals $X' - T'$ overconvergent along $X' \cap T'$, we have an isomorphism $\Db_{T'} (\RG^\dag_{X'} f^! \E [\delta_{X} - \delta_{X'}]) \simeq  \RG^\dag_{X'} f^! (\Db_{T'}\E) [\delta_{X} - \delta_{X'}]$. Then as in the proof of \cite[Proposition 13.4.5]{caro2021arithmetic}, one gets a split injection $\E \hookrightarrow f_+\RG_{X'}^\dag f^!\E$.
	\end{proof}
		
	\subsection{$\sp_+$ functor}
	We have constructed a functor $\sp_{X \hookrightarrow \Pf, T} : \MIC^\dag(X-T, X, \Pf / K) \simto \MIC^\ddag(X, \Pf, T / K)$ for a formal $\V[[t]]$-scheme of finite type locally having finite $p$-bases over $\V$, a closed subscheme $X$ of $P$ locally having finite $p$-bases over $k$ and a divisor $T$ of $\Pf$. In order to apply semistable reduction theorem to the study of overconvergent $F$-isocrystals and to the construction of six functor formalism of arithmetic $\D$-modules, we need to define a suitable functor $\sp_{X \hookrightarrow \Pf, T}$ to the case where $T$ is not necessarily a divisor of $P$.
	
	In what follows, we largely follow \cite{Caro2011} and \cite[\S16]{caro2025arith}. Although most of the proofs work similarly in our situation, we have to remark that de Jong's alteration theorem is used all over the place, which we should be careful. The $k[[t]]$-version of de Jong's alteration theorem \cite[Theorem 6.5]{de1996smoothness} is as follows : Given an integral, separated $k[[t]]$-scheme of finite type $X$ and a proper closed subscheme $Z$, there is a finite extension $k'[[t']]$ of $k[[t]]$ and a projective alteration $f : Y \to X$ over $k[[t]]$ such that $(Y, f^{-1}(Z)_\red)$ is a strictly semistable pair over $k'[[t']]$. Since $k((t))$ is not perfect, in general $a$ is not generically etale. However, we can choose some $n \geq 0$ so that $k'[[t']]$ is etale over $k[[t^{1/p^n}]]$. The base changes via $k[[t]] \to k[[t^{1/p^n}]]$ does not affect the category of isocrystals or arithmetic $\D$-modules (e.g. \cite[Corollary 1.1.4.8]{caro2021arithmetic}). With this subtle modification, we can adopt the proofs in our situation.
	
	We display two important theorems in order to follow the arguments.
	
	\begin{theorem}[Kedlaya contagiosity phenomena, {\cite[5.3.4]{Caro2011}}{\cite[16.1.8.1]{caro2025arith}}]
		Let $U \hookrightarrow X \hookrightarrow Y$ be open immersions with $X$ smooth over $k[[t]]$, $Y$ a $k[[t]]$-variety and $U$ dense in $Y$. Then the canonical functor
		$$ \Isoc(X, Y/K) \to \Isoc(X, X/K) \times_{\Isoc(U, X/K)} \Isoc(U, Y/K)$$
		is an equivalence of categories.
	\end{theorem}
	\begin{proof}
		If $U$ and $X$ are schemes over $k((t))$, the theorem for $\Ed_K$-valued overconvergent isocrystals is proved in \cite[Proposition 5.2.7]{liu2026semistable}. Its proof is based on the results of log $\nabla$-modules relative to $\Ed_K$ in \S3 of {\it loc. cit.} The same results for log $\nabla$-modules relative to $K$ also hold true as explained in \S7 of {\it loc. cit.}, and the same proof works.
	\end{proof}
	
	\begin{theorem}[Proper descent, {\cite[Proposition 7.3]{shiho2007relative}}] \label{thm:proper-descent}
		Let $(f^{(\bullet)}, g^{(\bullet)})(X^{(\bullet)}, Y^{(\bullet)}) \to (X, Y)$ be a strongly proper \v{C}ech hypercover of formal pairs over $k$. Then there is an equivalence of categories
		$$\Isoc(X, Y/K) \simto 2\mathchar`\-\varprojlim\left[ \Isoc(X^{(0)}, Y^{(0)}) \rightrightarrows \Isoc(X^{(1)}, Y^{(1)}/K) \mathrel{\substack{\longrightarrow \\[-0.6ex] \longrightarrow \\[-0.6ex] \longrightarrow}} \Isoc(X^{(2)}, Y^{(2)}/K) \right].$$
	\end{theorem}
	\begin{remark}
		A formal pair $(X, Y)$ over $k$ is a open immersion of $k$-schemes locally of formal finite type. A morphism $(f, g) : (X', Y') \to (X, Y)$ is a strongly proper cover if $g$ is proper surjective and  $g^{-1}(X) = X'$.
	\end{remark}
	\begin{proof}[{\rm (Proof of \ref{thm:proper-descent})}]
		As in the proof of \cite[Proposition 7.3]{shiho2007relative}, we may assume that there are compatible closed immersions $Y \hookrightarrow \Yf$ into a formally smooth formal $\V$-scheme and a closed immersion $Y^{(0)} \hookrightarrow \wh{\mathbb{P}}^n_\Yf$, which makes
		\[\xymatrix{
			Y^{(\bullet)} \ar@{^(->}[r] \ar[d]_{g^{(\bullet)}} & (\wh{\mathbb{P}}^n_\Yf)^{\times_{\Yf} (\bullet+1)} \ar[d] \\
			Y \ar@{^(->}[r] & \Yf
		}\]
		commute. By definition,
		$$\Isoc(X^{(\bullet)}, Y^{(\bullet)}/K) \simto  2\mathchar`\-\varprojlim\left[ \Coh(]X[_{(\wh{\mathbb{P}}^n_\Yf)^{\times_{\Yf} (\bullet+1)}}) \rightrightarrows \Coh(]X[_{((\wh{\mathbb{P}}^n_\Yf)^{\times_{\Yf} (\bullet+1)})^{\times 2}}) \mathrel{\substack{\longrightarrow \\[-0.6ex] \longrightarrow \\[-0.6ex] \longrightarrow}} \Coh(]X[_{((\wh{\mathbb{P}}^n_\Yf)^{\times_{\Yf} (\bullet+1)})^{\times 3}}) \right].$$
		Therefore, the 2-limit of this theorem is that of the following diagram:
		\[\xymatrix{
			\Coh(]X[_{(\wh{\mathbb{P}}^n_\Yf)}) \ar@<-0.5ex>[r]\ar@<0.5ex>[r] \ar@<-0.5ex>[d]\ar@<0.5ex>[d] & \Coh(]X[_{(\wh{\mathbb{P}}^n_\Yf)^{\times 2}}) \ar@<-0.75ex>[r]\ar@<0.75ex>[r] \ar[r] \ar@<-0.5ex>[d]\ar@<0.5ex>[d] & \Coh(]X[_{(\wh{\mathbb{P}}^n_\Yf)^{\times 3}}) \ar@<-0.5ex>[d]\ar@<0.5ex>[d] \\
			\Coh(]X[_{(\wh{\mathbb{P}}^n_\Yf)^{\times_{\Yf} 2}}) \ar@<-0.5ex>[r]\ar@<0.5ex>[r] \ar@<-0.75ex>[d]\ar@<0.75ex>[d]\ar[d] & \ar@<-0.75ex>[r]\ar@<0.75ex>[r] \ar[r] \Coh(]X[_{((\wh{\mathbb{P}}^n_\Yf)^{\times_{\Yf} 2})^{\times 2}}) \ar@<-0.75ex>[d]\ar@<0.75ex>[d]\ar[d] &  \Coh(]X[_{((\wh{\mathbb{P}}^n_\Yf)^{\times_{\Yf} 2})^{\times 3}}) \ar@<-0.75ex>[d]\ar@<0.75ex>[d]\ar[d]\\
			\Coh(]X[_{(\wh{\mathbb{P}}^n_\Yf)^{\times_{\Yf} 3}}) \ar@<-0.5ex>[r]\ar@<0.5ex>[r] & \ar@<-0.75ex>[r]\ar@<0.75ex>[r] \ar[r] \Coh(]X[_{((\wh{\mathbb{P}}^n_\Yf)^{\times_{\Yf} 3})^{\times 2}}) & \Coh(]X[_{((\wh{\mathbb{P}}^n_\Yf)^{\times_{\Yf} 3})^{\times 3}}).\\
		}\]
		Thus it suffices to prove that the 2-limit of $i$-th column is equivalent to $\Coh(]X[_{\Yf^{i}}).$ Let $\Pf^i = (\wh{\mathbb{P}}^n_\V)^{\times i}$ and $\Pf^i_{\Yf^i} = \Pf^i \times \Yf^i$. Then the $i$-th column is induced from the diagram
		$$ (\Pf^i_{\Yf^i})^{\times_{\Yf^i} 3} \mathrel{\substack{\longrightarrow \\[-0.6ex] \longrightarrow \\[-0.6ex] \longrightarrow}} (\Pf^i_{\Yf^i})^{\times_{\Yf^i} 2} \rightrightarrows \Pf^i_{\Yf^i}.$$
		Replacing $\Pf^i$ by $\wh{\mathbb{P}}^{ni}_\V$ does not change the tubes, so it suffices to consider the first column. Since $]X[_{(\wh{\mathbb{P}}^n_\Yf)^{\times_{\Yf} j}} = ]X[_{\Yf} \times ( ]\{0\}[_{\wh{\mathbb{P}}^n_\V})^j = ]X[_\Yf \times (A_K^n[0, 1))^j$, the descent is reduced to the local calculation, which is straightforward.
	\end{proof}
	
	\begin{_empty}
		Let $(Y, X, \Pf)$ be a formal frame with $\Pf$ of finite type over $\V[[t]]$, locally having finite $p$-bases over $\V$ and smooth over $\V[[t]]$ or over $\V$ in an open neighborhood of $Y$. Let $Z$ be a closed subscheme of $P$ such that $Y = X - Z$.
		We denote by $\MIC^{(\bullet)}(X, \Pf, Z / K)$ the full subcategory of $\underrightarrow{LD}^{b, \heartsuit}_{\Qb, {\rm poc}}(X, \Pf, Z / \V)$ consisting of objects $\E$ such that $\varinjlim \E|_{\Pf - Z} \in \Coh(\D^\dag_{\Pf - Z, \Qb})$ belongs to $\MIC^\ddag(Y, \Pf - Z, \empty / K)$. We get the same category if we use ${\rm povcoh}$ instead of ${\rm poc}$. 
		
		Recall that when $Z = T$ is a divisor of $P$, we have constructed the functor
		$$\sp_{X \hookrightarrow \Pf, T, +} : \MIC^\dag(Y, X, \Pf / K) \to \MIC^\ddag(X, \Pf, T / K).$$
		We can prove that $\MIC^\ddag(X, \Pf, T / K)$ is equivalent to  $\MIC^{(\bullet)}(X, \Pf, T / K)$ under $\varinjlim$ similarly as in \cite[Theorem 16.1.1.17]{caro2025arith}. 
		Following the method of \cite[Theorem 16.2.1.10]{caro2025arith}, one constructs an equivalence $\sp_{X \hookrightarrow \Pf, Z, +} : \MIC^\dag(Y, X, \Pf / K) \simto \MIC^{(\bullet)}(X, \Pf, Z / K)$. Note that $\MIC^{(\bullet)}(X, \Pf, Z / K)$ and $\sp_{X \hookrightarrow \Pf, Z, +}$ are independent of the choice of $X, Z$, only dependent on $Y = X - Z$ as in \cite[Proposition 16.2.7.4]{caro2025arith}.
	\end{_empty}
	\begin{proposition} \label{prop:sp-dual}
		There is a canonical functor $\Db_Z : \MIC^{(\bullet)}(X, \Pf, Z / K)^\op \to \MIC^{(\bullet)}(X, \Pf, Z / K)^\op$ such that 
		\begin{itemize}
			\item $\Db_Z \Db_Z \E \simto \E$ for any $\E \in \MIC^{(\bullet)}(X, \Pf, Z / K)$ and
			\item $\sp_{X \hookrightarrow \Pf, Z, +} (E^\vee) \simto \Db_Z \sp_{X \hookrightarrow \Pf, Z, +} (E)$ for any $E \in \MIC^\dag(X, \Pf, Z / K)$.
		\end{itemize}
	\end{proposition}
	\begin{proof}
		Similar to \cite[16.2.3]{caro2025arith}.
	\end{proof}
	
	\begin{_empty}
		Let $(Y', X', \Pf')$ be a formal frame with $\Pf'$ of finite type over $\V[[t]]$, locally having finite $p$-bases over $\V$ and smooth over $\V[[t]]$ or over $\V$ in an open neighborhood of $Y'$. Let $Z'$ be a closed subscheme of $P'$ such that $Y' = X' - Z'$. Let $\theta = (u, v, f) : (Y', X', \Pf') \to (Y, X, \Pf)$ be a morphism of formal frames over $\V$ and suppose that $Z' \supseteq f^{-1}(Z)$. As in \cite[16.2.4]{caro2025arith}, one defines the functor
		$$\theta^! := \RG^\dag_{Y'} f^{!}[-\delta_{X'}+\delta_X] : \MIC^{(\bullet)}(X, \Pf, Z / K) \to \MIC^{(\bullet)}(X', \Pf', Z' / K)$$
		and
		$$\theta^+ := \Db_{Z'} \theta^+ \Db_{Z}.$$
	\end{_empty}
	\begin{proposition} \label{prop:sp-pullback}
		We have the following:
		\begin{itemize}
			\item Let $(Y'', X'', \Pf'')$ be a formal frame with $\Pf''$ of finite type over $\V[[t]]$, locally having finite $p$-bases over $\V$ and smooth over $\V[[t]]$ or over $\V$ in an open neighborhood of $Y''$. Let $Z''$ be a closed subscheme of $P''$ such that $Y'' = X'' - Z''$. Let $\theta' = (u', v', f') : (Y'', X'', \Pf'') \to (Y', X', \Pf')$ be a morphism of formal frames over $\V$ and suppose that $Z'' \supseteq f'^{-1}(Z')$. Then we have canonical equivalences $\theta'^! \theta^! \simeq (\theta \theta')^!$ and $\theta'^+ \theta^+ \simeq (\theta \theta')^+$.
			\item $\theta^!$ and $\theta^+$ are canonically isomorphic.
			\item We have a canonical isomorphism $\theta^+ \sp_{X \hookrightarrow \Pf, Z, +} \simeq \sp_{X' \hookrightarrow \Pf', Z', +} \theta^*$.
		\end{itemize}
	\end{proposition}
	
	\begin{proposition} \label{prop:sp-tensor}
		We have the following:
		\begin{itemize}
			\item Let $\E, \F \in \MIC^{(\bullet)}(X, \Pf, Z / K)$. Then $\E \wh{\otimes}^{\mathbb{L}}_{\O_{\Pf}} \F [\delta_{Y} - \delta_P] \in \MIC^{(\bullet)}(X, \Pf, Z / K)$.
			\item Let $E, F \in \MIC^\dag(Y, X, \Pf / K)$. We have a canonical isomorphism $\sp_{X \hookrightarrow \Pf, Z, +}(E \otimes_{j^\dag_Y \O_{]X[_\Pf}} F) \simeq \sp_{X \hookrightarrow \Pf, Z, +}(E) \wh{\otimes}^{\mathbb{L}}_{\O_{\Pf}} \sp_{X \hookrightarrow \Pf, Z, +}(F) [\delta_{Y} - \delta_P]$.
		\end{itemize}
	\end{proposition}
	\begin{proof}
		Similar as in \cite[16.2.6]{caro2025arith}.
	\end{proof}
	
	\begin{proposition} \label{prop:bullet-indep-frame}
		Let $X', Z'$ be another closed subschemes of $P$ such that $Y = X' - Z'$. Then we have the equality $\MIC^{(\bullet)}(X, \Pf, Z / K) = \MIC^{(\bullet)}(X', \Pf, Z' / K)$. 
	\end{proposition}
	\begin{proof}
		Similar as in \cite[Proposition 16.2.7.4]{caro2025arith}.
	\end{proof}
	
	We will denote $\MIC^{(\bullet)}(Y, \Pf / K)$ instead of $\MIC^{(\bullet)}(X, \Pf, Z / K)$ for some closed subschemes $X, Z$ of $P$ such that $X - Z = Y$.
	
	\subsection{The category $\LD{\rm dev}$}
	\begin{definition}[{cf. \cite[16.3]{caro2025arith}}]
		$\Pf$ be a formal $\V[[t]]$-scheme of finite type locally having finite $p$-bases over $\V$ and $Y$ a locally closed subscheme of $P$.
		\begin{itemize}
			\item A $p$-smooth stratification of $Y$ is a descending chain of reduced closed subschemes $Y_{\red} = Y_0 \supseteq \cdots \supseteq Y_m = \emptyset$ such that $Y_i - Y_{i+1}$ locally have finite $p$-bases over $k$ for all $i = 0, \dots, m-1$.
			\item We denote by $\LD{isoc}(X, \Pf , Z / \V)$ the triangulated full subcategory of $\LD{poc}(X, \Pf, Z / \V)$ consisting of those objects $\E$ such that all $H^n_Z(\E)$ belongs to $\MIC^{(\bullet)}(X, \Pf, Z / K).$ Note that this is only dependent on $Y$ as in {\rm \cite[Notation 16.3.1.4]{caro2025arith}}, we will denote it by $\LD{isoc}(Y, \Pf / \V)$.
			\item We say that $\E \in \LD{qc}(\wh{\D}^{(\bullet)}_{\Pf})$ splits into overconvergent isocrystals over $Y$ if there is some $n \geq 0$ such that for $p : \Pf^{(n)} = \Pf \times_{\Spf \V[[t]]} \Spf \V[[t^{1/p^n}]] \to \Pf$ the projection, $\E^{(n)} = p^! \E \simeq p^+ \E$ and $Y^{(n)} = Y \times_{\Spec k[[t]]} \Spec k[[t^{1/p^n}]]$, there is a $p$-smooth stratification $Y_0 \supseteq \cdots \supseteq Y_m = \emptyset$ of $Y^{(n)}$ such that $\RG^\dag_{Y_i} \E^{(n)}$ belongs to $\LD{isoc}(Y_i, \Pf^{(n)} / \V).$
			\item When $\E \in \LD{qc}(\wh{\D}^{(\bullet)}_{\Pf})$ splits into overconvergent isocrystals over $P$, we say that $\E$ is devissable in isocrystals. The full subcategory of $\LD{qc}(\wh{\D}^{(\bullet)}_{\Pf})$ consisting of objects devissable in isocrystals is denoted by $\LD{dev}(\wh{\D}^{(\bullet)}_{\Pf})$. We denote by $\LD{dev}(Y, \Pf / \V)$ the full subcategory of $\LD{dev}(\wh{\D}^{(\bullet)}_{\Pf})$ consisting of objects $\E$ such that there is an isomorphism $\RG^\dag_{Y}\E \simeq \E$.
		\end{itemize}
	\end{definition}
	
	\begin{proposition} \label{prop:povcoh-dev}
		Let $(Y, X, \Pf)$ be a formal frame with $\Pf$ of finite type over $\V[[t]]$, locally having finite $p$-bases over $\V$ and smooth over $\V[[t]]$ or over $\V$ in an open neighborhood of $Y$. Then $\LD{povcoh}(X, \Pf, Z / \V) \subseteq \LD{dev}(Y, \Pf / \V)$.
	\end{proposition}
	\begin{proof}
		When $\Pf$ is smooth over $\V$ in an open neighborhood of $Y$ and (hence) $Y$ is smooth over $k$, it is proved in the similar way as in \cite[Theorem 16.3.1.8]{caro2025arith}. When $\Pf$ is smooth over $\V[[t]]$ in an open neighborhood of $Y$ and (hence) $Y$ is smooth over $k[[t]]$, \cite[Proposition 14.1.11]{caro2021arithmetic} ensures devissability.
	\end{proof}
	
	\begin{proposition} \label{prop:dev-tensor}
		Let $(Y, X, \Pf)$ be a formal frame with $\Pf$ of finite type over $\V[[t]]$, locally having finite $p$-bases over $\V$ and smooth over $\V[[t]]$ or over $\V$ in an open neighborhood of $Y$. Let $Z$ be a closed subscheme of $P$ such that $Y = X - Z$. Then we have the functor
		$$- \wh{\otimes}^{\mathbb{L}}_{\O_\Pf} - : \LD{dev}(Y, \Pf / \V) \times \LD{dev}(Y, \Pf / \V) \to \LD{dev}(Y, \Pf / \V)$$
	\end{proposition}
	\begin{proof}
		Similar to \cite[Theorem 16.3.2.2]{caro2025arith}.
	\end{proof}
	
	\section{A comparison theorem between log rigid cohomology and rigid cohomology}
	\subsection{The proof of the comparison theorem}
	\begin{definition} \label{def:ssp}
		A strictly semistable pair over $\V[[t]]$ is a pair $(\Xf, \Zf)$ of formal $\V[[t]]$-schemes $\Xf, \Zf$ with $\Zf$ a closed formal subscheme of $\Xf$ such that Zariski locally on $\Xf$ there is an etale morphism of the form $\Xf \to \Spf \V[[t]]\langle t_1, \dots, t_n \rangle/(t-t_1 \cdots t_r)$ and $\Zf$ is the pullback of
		$V(t_1 \cdots t_m)$ for some $r \leq m \leq d$. Note that $\Zf$ is a strict normal crossings divisor of $\Xf$ relative to $\V$ and it defines a log structure on $\Xf$.
	\end{definition}
	\begin{_empty}
		Let $(\Tf, \Tf')$ be a strictly semistable pair over $\V[[t]]$ (resp. a smooth pair over $\V$), $g : \Xf \to \Tf$ a smooth morphism (resp. a morphism locally having finite $p$-bases), $\Zf \subseteq \Xf$ a closed formal subscheme such that $(\Xf, \Zf)$ is a strictly semistable pair over $\V[[t]]$ and that $\Xf' := g^{-1}\Tf' \subset \Zf$. We denote by $\Xf^\sharp = (\Xf, \Zf)$ and $\Tf^\sharp = (\Tf, \Tf')$ the log formal schemes over $\Bf^\sharp := (\Spf \V[[t]], V(t))$. By assumption, $g : \Xf \to \Tf$ induces the morphism of log formal schemes $g^\sharp : \Xf^\sharp \to \Tf^\sharp$. Let $\Zf' \subseteq \Zf$ be the union of irreducible components of $\Zf$ not contained in $\Xf'$. Let $D$ be a divisor of $X$, $U = X-D$ its complement, and let $U^\sharp, X^\sharp$ be the log schemes induced by $\Xf^\sharp$. Recall that an overconvergent log isocrystal $E$ over $(U^\sharp, X^\sharp, \Xf^\sharp)/\Bf^\sharp$ is a locally free $j^\dag_{U}\O_{\Xf_K}$-module $E$ of finite rank with an integrable log connection $\nabla : E \to E \otimes \Omega^1_{\Xf^\sharp_K / \Bf^\sharp}$ with the overconvergent condition.
		
		For an irreducible component $Z_i$ of $Z'$, by \cite[Lemma 7.1.1]{liu2026semistable} locally on $\Xf$ there is a canonical isomorphism $]Z_i[_{\Xf} \simeq \Zf_{i, K} \times_K A^1_K[0, 1)$ whose projection to the annulus is given by the local section cutting out $\Zf_i$ in $\Xf$. Let $E$ be an overconvergent log isocrystal on $(U^\sharp, X^\sharp, \Xf^\sharp)$ coming from a log $\nabla$-module over some open neighborhood $\Us$ of $]U[_\Xf$ . When $Z_i$ is not contained in $D$, $\E$ restricts to a log $\nabla$-module over $\Vs \times A^1_K[0, 1)$ for the corresponding open neighborhood $\Vs$ of $]Z_i \cap U[_{\Zf_i}.$ There is a unique generic point $\Spa L_i$ of $\Zf_{i, K}$ and $\E$ induces a log $\nabla$-module $E_i$ over $A^1_L[0, 1)$. This is independent of the choice of $\E$. The eigenvalues of the residue map of $E_i$ are called exponents of $E$ along $Z_i$. A priori, they belong to $L_i$; by the overconvergent condition on $E$, they belong to $\mathbb{Z}_p$. We pose the following assumption on the exponents:
	\end{_empty}
	
	\begin{hypothesis}
		Let $Z_1, \dots, Z_r$ be irreducible components of $Z'$ not contained in $D$.
		\begin{itemize}
			\item (Exponents condition) The exponents of $E$ along $Z_i$ are all equal to a single integer $e_i \in \mathbb{Z}.$
			\item (Strong exponents condition) Furthermore $e_i \leq 0$ for all $i$.
		\end{itemize}
		
	\end{hypothesis}
	
	\begin{_empty}
		If $Z_i \not\subseteq D$ and $\I_i$ is the sheaf of ideals of $Z_i$ in $X$, let $\I_{i, K} := \sp^* \I_{i}$; then $j^\dag_{U}\I_{i, K}$ is an overconvergent log isocrystal over $(U^\sharp, X^\sharp, \Xf^\sharp)/\Bf^\sharp$. We define $E(\sum m_i \Zf_i)$ to be the overconvergent log isocrystal corresponding to $E \otimes \bigotimes_i j^\dag_{U}\I_{i, K}^{\otimes -m_i}.$ The exponents of $E(\sum m_i \Zf_i)$ along $Z_j$ are those of $E$ along $Z_j$ minus $m_i$. For $Y = X - Z$, there is a canonical commutative diagram
		\begin{align}
			\label{diag:twist}
			\xymatrix{
				E \ar@{^{(}->}[r] \ar[dr] & E(\sum_i m_i \Zf_i) \ar[d] \\
				& j^\dag_{Y \cap Z} \E.
			}
		\end{align}
		
		Under these settings, we are able to prove following theorems:
	\end{_empty}
	
	\begin{theorem} \label{thm:compA}
		For any integers $m_i \geq \max(0, e_i)$, the diagram (\ref{diag:twist}) induces the canonical isomorphism
		$$\mathbb{R}g_{K*} \ul\Gamma^\dag_{]Z[}(j^\dag_U\Omega^\bullet_{\Xf_K^\sharp/\Tf^\sharp_K} \otimes_{j^\dag_U \O_{\Xf_K}} E) \xrightarrow{\sim} \allowbreak \mathbb{R}g_{K*} \Cone[j^\dag_U\Omega^\bullet_{\Xf_K^\sharp/\Tf^\sharp_K} \otimes_{j^\dag_U \O_{\Xf_K}} E \to \allowbreak j^\dag_U\Omega^\bullet_{\Xf_K^\sharp/\Tf^\sharp_K} \otimes_{j^\dag_U \O_{\Xf_K}} E(\sum m_i \Zf_i)][-1].$$
	\end{theorem}
	\begin{proof}
		As in the proof of \cite[Theorem 1.1.4]{caro2012overholonomicity}, it suffices to prove that if $m_i = 0 \geq e_i$ for all $i$, then
		$$\mathbb{R}g_{K*} (\Omega^\bullet_{\Xf_K^\sharp/\Tf^\sharp_K} \otimes_{\O_{\Xf_K}} E) \to \mathbb{R}g_{K*} j^\dag_{Y} (\Omega^\bullet_{\Xf_K^\sharp/\Tf^\sharp_K} \otimes_{\O_{\Xf_K}} E)$$
		is an isomorphism. Also from {\it loc. cit.}, the assertion is local on $\Xf$ and $\Tf$, we assume that $\Xf$, $\Tf$ are affine, there is a $p$-basis $z_1, \dots, z_d$ of $\Xf$ over $\Tf$ and $\Zf$ is given by $z_1 \cdots z_{r'} = 0$. Let $\Zf_i$ be the irreducible component of $\Zf$ cut out by $z_i$ and suppose that $Z_1, \dots, Z_r$ are those not contained in $D$.
		
		Put $V^\eta = \{ x \in \Xf_K : |z_i|_{[x]} > \eta \textrm{ for } i = 1, \dots, r \}$. Then $\{ V^\eta \}_\eta$ is a cofinal system of open neighborhoods of $]Y \cap U[_{\Xf}$ in $\Xf_K.$ Hence for any rational open subset $\Vs$ of $\Tf_K$ and a complex $\F$ of coherent $j^\dag_{U} \O_{\Xf_K}$-modules, we have
		$$\Rb\Gamma(\Vs, \Rb g_{K*} j^\dag_{Y} \F) \simeq \varinjlim_{\eta} \Rb\Gamma(g_K^{-1}(\Vs) \cap V^\eta_{}, \F).$$
		Hence it suffices to fix $\eta$ close enough to $1$ and prove that
		$$\Rb\Gamma(g_K^{-1}(\Vs), \Omega^\bullet_{\Xf_K^\sharp/\Tf^\sharp_K} \otimes_{\O_{\Xf_K}} E) \to \Rb\Gamma(g_K^{-1}(\Vs) \cap V^\eta_{}, \Omega^\bullet_{\Xf_K^\sharp/\Tf^\sharp_K} \otimes_{\O_{\Xf_K}} E)$$
		is an isomorphism. We will show this in the style of \cite[Proposition 2.4.6]{shiho2002crystallineII}.
		For $I \subseteq \{ 1, \dots, r \},$ put
		$$ W^\eta_I := \{ x \in g_K^{-1}(\Vs) : |z_i|_{[x]} < 1 \textrm{ for } i \in I \textrm{ and } |z_i|_{[x]} > \eta \textrm{ for } i \in \{1, \dots, r\} - I \}, $$
		$$ U^\eta_I := \{ x \in g_K^{-1}(\Vs) : |z_i|_{[x]} = 0 \textrm{ for } i \in I \textrm{ and } |z_i|_{[x]} > \eta \textrm{ for } i \in \{1, \dots, r\} - I \}, $$
		and
		$$V^\eta_I := W^\eta_I \cap V^\eta.$$
		For $I_0, \dots, I_p \subseteq \{ 1, \dots, r \}$, we have
		$$ W^\eta_{I_0} \cap \cdots \cap W^\eta_{I_p} = \{ x \in g_K^{-1}(\Vs) : |z_i|_{[x]} < 1 \textrm{ for } i \in I_0 \cup \cdots \cup I_p \textrm{ and } |z_i|_{[x]} > \eta \textrm{ for } i \in \{1, \dots, r\} - (I_0 \cap \cdots \cap I_p) \}. $$
		Put
		$$ U'^\eta_{I_0, \dots, I_p} := \{ x \in g_K^{-1}(\Vs) : |z_i|_{[x]} = 0 \textrm{ for } i \in I_0 \cup \cdots \cup I_p \textrm{ and } |z_i|_{[x]} > \eta \textrm{ for } i \in \{1, \dots, r\} - (I_0 \cap \cdots \cap I_p) \}.$$
		Then we have the commutative diagram
		\[
		\xymatrix{
			& W^\eta_{I_0 \cup \cdots \cup I_p} \ar[rr]^{\sim} &  & U^\eta_{I_0 \cup \dots \cup I_p} \times A^{|I_0 \cup \cdots \cup I_p|}[0, 1) \\	
			W^\eta_{I_0} \cap \cdots \cap W^\eta_{I_p} \ar@{^(->}[ru] \ar[rr]^{\sim}&  & U'^\eta_{I_0, \dots, I_p} \times A^{|I_0 \cup \cdots \cup I_p|}[0, 1) \ar@{^(->}[ru] &  \\
			& V^\eta_{I_0 \cup \cdots \cup I_p} \ar@{^(-->}[uu] \ar@{-->}[rr]^{\sim}&  & U^\eta_{I_0 \cup \dots \cup I_p} \times A^{|I_0 \cup \cdots \cup I_p|}(\eta, 1) \ar@{^(->}[uu]\\
			V^\eta_{I_0} \cap \cdots \cap V^\eta_{I_p} \ar@{^(->}[uu] \ar@{^(->}[ru] \ar[rr]^{\sim} &  & U'^\eta_{I_0, \dots, I_p} \times A^{|I_0 \cup \cdots \cup I_p|}(\eta, 1) \ar@{^(->}[uu] \ar@{^(->}[ru] &  \\
		}
		\]
		whose projections to the annulus is given by those $z_i$ in $I_0, \dots, I_p$. $g_K^{-1}(\Vs)$ is covered by $\{ W^\eta_I \}_I$ and $g_K^{-1}(\Vs) \cap V^\eta$ is covered by $\{ V^\eta_I \}_I$. Using the spectral sequences induced by these coverings, it suffices to prove that
		$$H^*(W^\eta_{I_0} \cap \cdots \cap W^\eta_{I_p}, \Omega^\bullet_{\Xf_K^\sharp/\Tf^\sharp_K} \otimes_{\O_{\Xf_K}} E) \to H^*(V^\eta_{I_0} \cap \cdots \cap V^\eta_{I_p}, \Omega^\bullet_{\Xf_K^\sharp/\Tf^\sharp_K} \otimes_{\O_{\Xf_K}} E)$$
		are isomorphisms. Let $s = |I_0 \cup \cdots \cup I_p|$. The case where $s = 0$ is trivial, so we assume that $s > 0$. Put $U' = U'^\eta_{I_0, \dots, I_p}$. We have to prove that
		$$ H^*(U' \times A^s[0, 1), \Omega^\bullet_{U'^\sharp \times A^s[0, 1)^\sharp/\Tf^\sharp_K} \otimes_{\O_{\Xf_K}} E) \to H^*(U' \times A^s(\eta, 1), \Omega^\bullet_{U'^\sharp \times A^s[0, 1)^\sharp/\Tf^\sharp_K} \otimes_{\O_{\Xf_K}} E) $$
		is an isomorphism. Note that hiher cohomologies of coherent $j^\dag_U \O$-modules over $U'$, $U' \otimes A^s[0, 1)$ or $U' \otimes A^s(\eta, 1)$ vanish, since the spaces are quasi-Stein, see \cite[Theorem A.4]{liu2026semistable}. Thus using Katz-Oda spectral sequences
		\begin{align*}
			 E_2^{p, q} & = H^p(\Gamma(U', \Omega^\bullet_{U'^\sharp/\Tf^\sharp_K}) \otimes H^q(U' \times A^s[0, 1), \Omega^\bullet_{U'^\sharp \times A^s[0, 1)^\sharp / U'^\sharp} \otimes E))\\
			 & \Rightarrow H^{p+q}(U' \times A^s[0, 1), \Omega^\bullet_{U'^\sharp \times A^s[0, 1)^\sharp/\Tf^\sharp_K} \otimes E)
		\end{align*}
		and
		\begin{align*}
			'E_2^{p, q} & = H^p(\Gamma(U', \Omega^\bullet_{U'^\sharp/\Tf^\sharp_K}) \otimes H^q(U' \times A^s(\eta, 1), \Omega^\bullet_{U'^\sharp \times A^s(\eta, 1) / U'^\sharp} \otimes E))\\
			& \Rightarrow H^{p+q}(U' \times A^s(\eta, 1), \Omega^\bullet_{U'^\sharp \times A^s(\eta, 1)^\sharp/\Tf^\sharp_K} \otimes E),
		\end{align*}
		it suffices to prove that
		$$ H^*(U' \times A^s[0, 1), \Omega^\bullet_{U'^\sharp \times A^s[0, 1)^\sharp / U'^\sharp} \otimes E) \to  H^*(U' \times A^s(\eta, 1), \Omega^\bullet_{U'^\sharp \times A^s(\eta, 1) / U'^\sharp} \otimes E)$$
		are isomorphisms. By exponents condition, a suitable twist of $E$ have nilpotent residues. By \cite[Proposition 3.5.3]{liu2026semistable}, it is unipotent, so $E$ is an iterated extension of log $\nabla$-modules $F \otimes \O \prod_i z_i^{e_i}$ with $F$ a constant log $\nabla$-module. By deivssage and taking an open cover of $U'$, we may assume $E = \O \prod_i z_i^{e_i}$. Now using $e_i \leq 0$, the assertion is verified by an easy calculation of the cohomologies.
	\end{proof}
	
	\begin{theorem} \label{thm:compB}
		Let $\Zf'_1$ be the union of all irreducible components of $\Zf$ other than a fixed one $\Zf_1$ whose special fiber $Z_1$ is not contained in $D$. Endow $\Zf_1$ with the log structure of $\Zf_1 \cap \Zf'_1$. Let $m_1 \geq \max(0, e_1)$ and let $\nabla_1 : E(m_1 \Zf_1)/E \to \Omega^1_{\Xf^\sharp_K/\Zf^\sharp_{1,K}} \otimes E(m_1 \Zf_1)/E$ be the log connection induced from $\nabla$.
		Suppose that $g_1 : \Xf \to \Zf_1$ over $\Tf$ is a morphism locally with finite $p$-basis which is a retraction of the closed immersion $\Zf_1 \hookrightarrow \Xf$ such that $g_1^{-1}(\Zf_1 \cap \Zf') = \Zf'_1.$ Then we have the natural isomorphism
		$$\mathbb{R}g_{1,K*} \ul\Gamma^\dag_{]Z_1[}(\Omega^\bullet_{\Xf_K^\sharp/\Zf_{1, K}^\sharp} \otimes E) \xrightarrow{\sim} [g_{1K*}(E(m_1 \Zf_1)/E) \xrightarrow{g_{1K*}\nabla_1} g_{1K*}(\Omega^1_{\Xf^\sharp_K/\Zf^\sharp_{1,K}} \otimes E(m_1 \Zf_1)/E)][-1].$$
		Furthermore, the Gauss-Manin connections makes $g_{1K*}(E(m_1 \Zf_1)/E)$ and $g_{1K*}(\Omega^1_{\Xf^\sharp_K/\Zf^\sharp_{1,K}}E(m_1 \Zf_1)/E)$ overconvergent log isocrystals over $((Z_1\cap U)^\sharp, Z_1^\sharp, \Zf_1^\sharp)/\Bf^\sharp$ satisfying the exponents condition.
	\end{theorem}
	\begin{proof}
		By Theorem \ref{thm:compA}, we have the natural isomorphism
		$$\mathbb{R}g_{1K*} \ul\Gamma^\dag_{]Z_1[}(\Omega^\bullet_{\Xf_K^\sharp/\Zf^\sharp_{1,K}} \otimes E) \xrightarrow{\sim} \mathbb{R}g_{1, K*} \Cone[\Omega^\bullet_{\Xf_K^\sharp/\Zf^\sharp_{1, K}} \otimes E \to \Omega^\bullet_{\Xf_K^\sharp/\Zf^\sharp_{1, K}} \otimes E(m_1 \Zf_1)][-1].$$
		Since
		$$ \Cone[\Omega^\bullet_{\Xf_K^\sharp/\Zf^\sharp_{1, K}} \otimes E \to \Omega^\bullet_{\Xf_K^\sharp/\Zf^\sharp_{1, K}} \otimes E(m_1 \Zf_1)] \simeq [E(m_1 \Zf_1)/E \to \Omega^1_{\Xf_K^\sharp/\Zf^\sharp_{1, K}} \otimes E(m_1 \Zf_1)/E], $$
		we have to prove that $E(m_1 \Zf_1)/E$ and $\Omega^1_{\Xf_K^\sharp/\Zf^\sharp_{1, K}} \otimes E(m_1 \Zf_1)/E$ are $g_{1, K*}$-acyclic and that $g_{1, K*}(E(m_1 \Zf_1)/E)$ is locally free. As in $0^\circ$ of the proof of \cite[Theorem 1.1.4]{caro2012overholonomicity}, $i_{1, K*}i_{1, K}^* E(\Zf_1) = E(\Zf_1)/E$ and hence $\Rb g_{1, K*} E(\Zf_1)/E = \Rb g_{1, K*} i_{1, K*}i_{1, K}^* E(\Zf_1) = i_{1,K}^* E(\Zf_1)$ is locally free. By induction and devissage, we get the assertion. 
	\end{proof}
	
	\subsection{Interpretation of the comparison theorem with arithmetic $\D$-modules}
	In order for connect overconvergent log isocrystals with arithmetic $\D$-modules, we need a logarithmic version of arithmetic $\D$-module $\D^\dag_{\Xf^\sharp/\V}(^\dag D)_\Qb$ for a log formal $\V[[t]]$-scheme $\Xf^\sharp$ induced by a strictly semistable pair $(\Xf, \Zf)$ over $\V[[t]]$ and a divisor $D$ of $X$. As for the classical case, which means that the case where $(\Xf, \Zf)$ is a smooth pair over $\V$, the theory is fully developed. Although the strictly semistable case is only partly developed in \cite{caro2022logarithmic}, we still use the same results as in the classical case. As remarks for doubtful readers, in this paper we only use log arithmetic $\D$-modules $\D^\dag_{\Xf^\sharp/\V}(^\dag D)_\Qb$ with $\Xf$ affine and having a global finite $p$-basis, which allows us to reduce all results to a simple computation using global sections and coordinates. From this perspective, since local descriptions of the ring $\D^\dag_{\Xf^\sharp/\V}(^\dag D)_\Qb$ and various functors such as extraordinary pullbacks, pushforwards and local cohomological functors are just the same as the classical case, one should easily prove statements of log arithmetic $\D$-modules of strictly semistable pairs as well. Henceforth, we will use the notations and propositions about log arithmetic $\D$-modules from \cite{caro2012overholonomicity}.
		
	\begin{proposition}[cf. {\cite[Proposition 1.3.1]{caro2012overholonomicity}}] \label{prop:isoc-Dmod-equiv}
		Let $(\Xf, \Zf)$ be a strictly semistable pair (resp. a smooth pair) with $\Xf$ affine and having a finite $p$-basis (resp. an etale coordinate), let $D$ be a divisor of $X$ and let $U$ be its complement. Let $\Xf^\sharp$ be the induced formal log scheme and let $X^\sharp$ and $U^\sharp$ be the induced log schemes. Then the functors $\sp_*$ and $\sp^*$ induce quasi-inverse equivalences between the category of coherent $\D^\dag_{\Xf^\sharp/\V}(^\dag D)_\Qb$-modules locally projective over $\O_{\Xf, \Qb}$ of finite rank and the category of overconvergent log isocrystals over $(U^\sharp, X^\sharp, \Xf^\sharp)/K.$
	\end{proposition}
	\begin{proof}
		The case of smooth pairs is proved in \cite[Proposition 1.3.1]{caro2012overholonomicity}. The strictly semistable case is proved in the similar way as the classical case by using the coordinate.
	\end{proof}
		
	\begin{proposition}[cf. {\cite[Lemma 1.3.5]{caro2012overholonomicity}}] \label{prop:Rsp}
		Let $(\Xf, \Zf)$ be a strictly semistable pair (resp. a smooth pair) with $\Xf$ affine and having a finite $p$-basis (resp. an etale coordinate), let $D$ be a divisor of $X$ and let $U$ be its complement. Let $E$ be an overconvergent log isocrystal over $(U^\sharp, X^\sharp, \Xf^\sharp)/K$ and let $\E = \sp_*E$ be the corresponding coherent $\D^\dag_{\Xf^\sharp/\V}(^\dag D)_\Qb$-modules locally projective over $\O_{\Xf, \Qb}$ of finite rank. For a divisor $D' \subseteq D$ and its complement $U'$, we have the following:
		\begin{enumerate}
			\item $\E(^\dag D') = \D^\dag_{\Xf^\sharp/\V}(^\dag D')_\Qb \otimes_{\D^\dag_{\Xf^\sharp/\V}(^\dag D)_\Qb} \E \simeq \sp_* (j^\dag_{U'} E).$
			\item $\RG^\dag_{D'}(\E) \simeq \Rb \sp_* (\ul\Gamma^\dag_{]D'[_{\Xf'}} E).$
		\end{enumerate}
	\end{proposition}
	
	\begin{theorem}[cf. {\cite[Theorem 1.3.6]{caro2012overholonomicity}}] \label{thm:comp-arithA}
		Let $(\Xf, \Zf)$, $(\Tf, \Tf')$ and $g : \Xf \to \Tf$ be one of the following:
		\begin{itemize}
			\item $(\Xf, \Zf)$ is a strictly semistable pair over $\V[[t]]$, $(\Tf, \Tf')$ is a strictly semistable pair over $\V[[t]]$ and $g$ is a smooth morphism.
			\item $(\Xf, \Zf)$ is a strictly semistable pair over $\V[[t]]$, $(\Tf, \Tf')$ is a smooth pair over $\V$ and $g$ is a morphism having finite $p$-basis.
			\item $(\Xf, \Zf)$ is a is a smooth pair over $\V$, $(\Tf, \Tf')$ is a smooth pair over $\V$ and $g$ is a smooth morphism.
		\end{itemize}
		Moreover suppose that $\Xf' := g^{-1}\Tf' \subset \Zf$ so that $g$ induces the morphism of log formal schemes $g^\sharp : \Xf^\sharp = (\Xf, \Zf) \to \Tf^\sharp = (\Tf, \Tf')$. Let $D$ be a divisor of $X$ and let $\E$ be a coherent $\D^\dag_{\Xf^\sharp/\V}(^\dag D)_\Qb$-modules locally projective over $\O_{\Xf, \Qb}$ of finite rank satisfying the strong exponents condition. Then the canonical map
		$$ \Rb g_* (\Omega^\bullet_{\Xf^\sharp/\Tf^\sharp, \Qb} \otimes_{\O_{\Xf, \Qb}} \E) \to \Rb g_* (\Omega^\bullet_{\Xf/\Tf, \Qb} \otimes_{\O_{\Xf, \Qb}} \E(^\dag Z))$$
		is an isomorphism. Equivalently, for $g^\sharp : \Xf^\sharp \to \Tf^\sharp$ the morphism of log formal schemes from $\Xf^\sharp = (\Xf, \Zf)$ and $\Tf^\sharp = (\Tf, \Tf')$, we have $g^\sharp_+ \RG^\dag_Z \E = 0.$
	\end{theorem}
	\begin{proof}
		All cases can be proved in the same way as \cite[Theorem 1.3.6]{caro2012overholonomicity} by applying $\Rb \sp_*$ to Theorem \ref{thm:compA} and use Proposition \ref{prop:Rsp}. The last statement is from \cite[Remark 1.3.7]{caro2012overholonomicity}.
	\end{proof}
	
	\begin{_empty} \label{emp:rho}
		As in \cite[1.3.8]{caro2012overholonomicity}, for $u : \Xf^\sharp \to \Xf$ the morphism forgetting the log structure, we have the functor
		$$u_{D, +} \E = \D^\dag_{\Xf \leftarrow \Xf^\sharp}(^\dag D)_{ \Qb} \Lotimes_{\D^\dag_{\Xf^\sharp}(^\dag D)_\Qb} \E.$$
		As in \cite[6.2.1]{caro2006}, $\D^\dag_{\Xf \leftarrow \Xf^\sharp}(^\dag D)_{ \Qb} \simeq \D^\dag_{\Xf}(^\dag D)_\Qb \otimes_{\O_{\Xf}} \O_{\Xf}(\Zf)$ and one has
		\begin{equation}
			\label{eq:u+}
			u_{D, +} \E \simeq \D^\dag_{\Xf}(^\dag D)_\Qb \otimes^{\mathbb{L}}_{\D^\dag_{\Xf^\sharp}(^\dag D)_\Qb} \E(\Zf).
		\end{equation}
		We have the canonical morphism
		$$\rho : u_{D, +} \E \to \D_{\Xf}(^\dag D \cup Z)_\Qb \Lotimes_{\D^\dag_{\Xf^\sharp}(^\dag D)_\Qb} \E = \E(^\dag Z)$$
		induced from the map $\D^\dag_{\Xf \leftarrow \Xf^\sharp}(^\dag D)_{ \Qb} \simto \D^\dag_{\Xf^\sharp}(^\dag D)_\Qb \otimes_{\O_\Xf} \O_{\Xf}(\Zf) \to \D_{\Xf}(^\dag D \cup Z)_\Qb.$ This morphism fits into the following exact triangle:
		$$ \RG^\dag_{Z} u_{D,+}\E \to u_{D, +}\E \xrightarrow{\rho} \E(^\dag Z) \to.$$
	\end{_empty}
	
	\begin{lemma} \label{lem:u+RGcom}
		Let $u : \Xf^\sharp \to \Xf$ and $D$ be as in \ref{emp:rho} and let $D'$ be a divisor of $X$. Then for $\E \in D(\D^\dag_{\Xf^\sharp}(^\dag D)_\Qb)$, there is a canonical isomorphism $u_{D, +} \RG^{\dag}_{D'} \E \simeq \RG^{\dag}_{D'} u_{D, +} \E$ in $D(\D^\dag_{\Xf}(^\dag D)_\Qb)$.
	\end{lemma}
	\begin{proof}
		See \cite[Lemma 1.3.11]{caro2012overholonomicity}.
	\end{proof}
	
	\begin{corollary} \label{cor:comp-arithB}
		In the situation of Theorem \ref{thm:comp-arithA}, let $u : \Xf^\sharp \to \Xf$ be the canonical morphism forgetting the log structure. Then $g_{D, +} (\rho) : g_{D, +}u_{D, +}\E \to g_{D, +} \E (^\dag Z)$ is an isomorphism and $g_{+} \RG_{Z}^\dag u_{D, +} \E = 0.$
	\end{corollary}
	\begin{proof}
		Similar to \cite[Theorem 1.3.13]{caro2012overholonomicity}.
	\end{proof}
	
	\begin{theorem} \label{thm:comp-arithB}
		Let $(\Xf, \Zf)$, $(\Tf, \Tf')$, $g : \Xf \to \Tf$ and $D$ be as in Theorem \ref{thm:comp-arithA}, let $\Zf'_1$ be the union of all irreducible components of $\Zf$ other than a fixed one $\Zf_1$ whose special fiber $Z_1$ is not contained in $D$. Endow $\Zf_1$ with the log structure via $\Zf \cap \Zf'_1$. Suppose that there is a retraction $g_1 : \Xf \to \Zf_1$ of the closed immersion $i_1 : \Zf_1 \hookrightarrow \Xf$ such that $g_1^{-1}(\Zf_1 \cap \Zf') = \Zf'.$ Let $g_1^\sharp : \Xf^\sharp \to \Zf_1^\sharp$ be the induced morphism of log formal schemes. For a coherent $\D^\dag_{\Xf^\sharp}(^\dag D)_\Qb$-module $\E$ locally projective of finite rank over $\O_{\Xf}(^\dag D)_\Qb$ satisfying the exponents condition, the complex
		$$g^\sharp_{1, +}\RG^\dag_{Z_1}\E \simeq \Cone[g_{1,+}^\sharp\E \to g_{1, +}^\sharp\E(^\dag Z_1)][-1]$$
		is isomorphic to a complex of coherent $\D^\dag_{\Zf_1^\sharp}(^\dag D\cap Z_1)_\Qb$-modules $\E$ locally projective of finite rank over $\O_{\Zf_1}(^\dag D \cap Z_1)_\Qb$ satisfying the exponents condition.
	\end{theorem}
	\begin{proof}
		Similar to \cite[Theorem 1.3.15]{caro2012overholonomicity}.
	\end{proof}
	
	\section{Application to the study of overconvergent $F$-isocrystals}
	\subsection{A comparison between log de Rham and de Rham complexes}
	\begin{_empty} \label{p-smooth-retraction}
		Let $\phi : \Xf \to \Spf \V[[t]]\langle t_1, \dots, t_n \rangle / (t - t_1 \cdots t_r)$ be an affine etale morphism, $\Zf_i = \phi^{-1}(V(t_i))$ for $i = 1, \dots, m$, where $r \leq m \leq n$ and let $R = \V[[t]]\langle t_1, \dots, t_r \rangle / (t - t_1 \cdots t_r)$. For $1 \leq j_1 < \cdots < j_k \leq m$, put $\Zf_{j_1, \dots, j_k} = \Zf_{j_1} \cap \cdots \cap \Zf_{j_k}$. If $r < j_1$, then we have the following commutative diagram:
		\begin{align*}
			\xymatrix{
				\Zf_{j_1, \dots, j_k} \ar@{_(->}[rd]^-{\iota'} \ar@{_(->}@/_10pt/[rdd]_{\iota} \ar@{_(->}@/^10pt/[rrd]^{(0, \id)}& & & \\
				& \Xf' \ar[d]_{f} \ar[r]^-{h} & \widehat{\mathbb{A}}^{k}_R \times_R \Zf_{j_1, \dots, j_k} \ar[r]^-{\mathrm{pr}} \ar[d]^{\id \times \phi'}& \Zf_{j_1, \dots, j_k}\\
				& \Xf \ar[r]^-{\phi}& \widehat{\mathbb{A}}^{k}_R \times_R \widehat{\mathbb{A}}^{n-r-k}_R, & 
			}
		\end{align*}
		where the square is cartesian, $\iota : \Zf_{j_1, \dots, j_k} \hookrightarrow \Xf$ is the closed immersion, $0 : \Zf_{j_1, \dots, j_k} \to \mathbb{A}^k_R$ is the zero section and $\phi' : \Zf_{j_1, \dots, j_k} \to \wh{\mathbb{A}}^{n-r-k}_R$ is the base change of $\phi$. 
		Note that $f$ and $h$ are etale. As in \cite[3.1.2]{abe2018theory}, by replaing $\Xf'$ by an affine open subset, we may assume that $f^{-1}(\Zf_{j_1} \cap \cdots \cap \Zf_{j_k}) = \iota'(\Zf_{j_1} \cap \cdots \cap \Zf_{j_k}).$ Let $g = \mathrm{pr} \circ h : \Xf' \to \Zf_{j_1} \cap \cdots \cap \Zf_{j_k}.$ Then $g$ is a smooth morphism and $g \circ \iota' = \id.$
		
		If $1 \leq j_1 < \cdots < j_l \leq r < j_{l+1} < \cdots < j_k \leq m$ with $l \geq 1$, then we have the following commutative diagram:
		\begin{align*}
			\xymatrix{
				\Zf_{j_1, \dots, j_k} \ar@{_(->}[rd]^-{\iota'} \ar@{_(->}@/_10pt/[rdd]_{\iota} \ar@{_(->}@/^10pt/[rrrd]^-{(0, \id)} & 	& &  & \\
				& \Xf'	 \ar[d]_{f} \ar[r]^-{h} & \Yf \ar@{^(->}[r]^-{i'} \ar[d]^{\psi} & \wh{\mathbb{A}}^{k}_{\V[[t]]} \times_\V \Zf_{j_1, \dots, j_k} \ar[d]^{\id \times \phi'} \ar[r]^-{\mathrm{pr}'} & \Zf_{j_1, \dots, j_k} \ar[d]^{\phi'} \\	
				& \Xf \ar[r]^-{\phi} & \Spf \V[[t]]\langle t_1, \dots, t_n \rangle / (t - t_1 \cdots t_r) \ar@{^(->}[r]^-{i} & \wh{\mathbb{A}}^{k}_{\V[[t]]} \times_\V \wh{\mathbb{A}}^{n-k}_{\V} \ar[r]^-{\mathrm{pr}} & \wh{\mathbb{A}}_\V^{n-k},
			}
		\end{align*}
		where the squares are all cartesian, $\iota$ and $\phi'$ are defined similarly as above, $0 : \Zf_{j_1, \dots, j_k} \to \wh{\mathbb{A}}^k_{\V[[t]]}$ is the composition of the zero section $\Zf_{j_1, \dots, j_k} \to \wh{\mathbb{A}}^k_{\V}$ with the closed immersion $\wh{\mathbb{A}}^k_{\V} \hookrightarrow \wh{\mathbb{A}}^k_{\V[[t]]}$. Again $f$ and $h$ are etale and we may assume $f^{-1}(\Zf_{j_1, \dots, j_k}) = \iota'(\Zf_{j_1, \dots, j_k}).$ Since the morphism $ \mathrm{pr} \circ i : \Spf \V[[t]]\langle t_1, \dots, t_n \rangle / (t - t_1 \cdots t_r) \to \wh{\mathbb{A}}_\V^{n-k}$ has a finite $p$-basis $t_{j_1}, \dots, t_{j_k}$, its base change $\mathrm{pr}' \circ i'$ also has a finite $p$-basis $\psi^*t_{j_1}, \dots, \psi^*t_{j_k}$. Let $g = \mathrm{pr}' \circ i' \circ h : \Xf' \to \Zf_{j_1, \dots, j_k}$. Then $g$ has a finite $p$-basis and $g \circ \iota' = \id$.
		
		In any case, there is an affine etale morphism $f : \Xf' \to \Xf$ such that $f^{-1}(\Zf_{j_1, \dots, j_k}) \xrightarrow{\sim} \Zf_{j_1, \dots, j_k}$ is an isomorphism and a retraction $g : \Xf' \to \Zf_{j_1, \dots, j_k}$ with a finite $p$-basis of the closed immersion $\iota' : \Zf_{j_1, \dots, j_k} \hookrightarrow \Xf'$:
		
		\begin{align}
			\label{diag:ret}
			\xymatrix{
			\Zf_{j_1, \dots, j_k} \ar@{_(->}[rd]^-{\iota'} \ar@{_(->}@/_10pt/[rdd]_-{\iota} \ar@{=}@/^15pt/[rrd] &  &\\
			& \Xf' \ar[d]^-{f} \ar[r]^-{g} & \Zf_{j_1, \dots, j_k} \\
			& \Xf.
			}
		\end{align}
	\end{_empty}
	
	\begin{lemma} \label{lem:etale-pullback}
		Suppose that $(\Xf, \Zf)$ is an affine strictly semistable pair over $\V[[t]]$ or an affine smooth pair over $\V$, $D \subseteq X$ a divisor, $f : \Xf' \to \Xf$ an affine smooth morphism with an etale coordinate, $\Zf' := f^{-1}(\Zf)$, $D' = f^{-1}(D)$, $f^\sharp : \Xf'^\sharp  \to \Xf^\sharp$ the canonical morphism of log formal schemes and $\E$ a coherent $\D^\dag_{\Xf^\sharp}(^\dag D)_\Qb$-module locally projective of finite rank over $\O_{\Xf^\sharp}(^\dag D)_\Qb$ satisfying the exponent condition (resp. strong exponent condition). For $E$ the corresponding overconvergent log isocrystal, the pullback $f^{\sharp, *}E$ corresponds to $f_D^{\sharp, !}\E[\delta_X - \delta_{X'}]$, which is a coherent $\D^\dag_{\Xf'^\sharp}(^\dag D')_\Qb$-module locally projective of finite rank over $\O_{\Xf'^\sharp}(^\dag D')_\Qb$ satisfying the exponent condition (resp. strong exponent condition).
	\end{lemma}
	\begin{proof}
		We have to check that the pullback $f^{\sharp, *}$ of overconvergent log isocrystals and the extraordinary pullback $f^{\sharp, !}_D$ of log arithmetic $\D$-modules are compatible after shift under the equivalence Proposition \ref{prop:isoc-Dmod-equiv}. When $(\Xf, \Zf)$ is a smooth pair, this is a special case of \cite[Lemma 2.3.10]{caro2012overholonomicity}. When $(\Xf, \Zf)$ is a strictly semistable pair, the proof is similar (in fact it is much easier as $f$ has an etale coordiante).
	\end{proof}
	\begin{lemma} \label{lem:u+comm}
		Suppose that $(\Xf, \Zf)$ is an affine strictly semistable pair over $\V[[t]]$ or an affine smooth pair over $\V$, $D \subseteq X$ a divisor, $f : \Xf' \to \Xf$ an affine smooth morphism with an etale coordinate, $\Zf' := f^{-1}(\Zf)$, $D' = f^{-1}(D)$, and let $f^\sharp : \Xf'^\sharp  \to \Xf^\sharp$, $u : \Xf^\sharp \to \Xf$ and $u' : \Xf'^\sharp \to \Xf'$ be the canonical morphism of log formal schemes. Then for $\E \in D_{\rm coh}^b(\D_{\Xf^\sharp}({^\dag D})_\Qb)$, we have the canonical isomorphism in $D^b_{\rm coh}(\D_{\Xf'}({^\dag D'})_\Qb)$:
		$$ f^{!}_D u_{D, +} \E \xrightarrow{\sim} u'_{D', +} f^{\sharp, !}_D \E. $$
	\end{lemma}
	\begin{proof}
		When $(\Xf, \Zf)$ is a smooth pair, this is \cite[Theorem 13.2.3.7]{caro2025arith}. One can prove the strictly semistable case as well using the coordinate.
	\end{proof}
	
	\begin{lemma} \label{lem:log-D-dagger}
		Let $(\Xf, \wt{\Zf})$ be an affine strictly semistable pair over $\V[[t]]$ or an affine smooth pair over $\V$, let $\Zf_1, \dots, \Zf_n$ be the irreducible components of $\wt{\Zf}$, $\Zf = \Zf_1 \cup \cdots \cup \Zf_r$, $\Zf' = \Zf_{r+1} \cup \cdots \cup \Zf_n$. Let $\wt{\Xf}^\sharp = (\Xf, \wt{\Zf})$ and $\Xf^\sharp = (\Xf, \Zf)$ be log formal schemes. Then for a divisor $D$ of $X$, the canonical morphism $\D^\dag_{\wt{\Xf}^\sharp}(^\dag D\cup Z')_\Qb \to \D^\dag_{\Xf^\sharp}(^\dag D\cup Z')_\Qb$ is an isomorphism.
	\end{lemma}
	\begin{proof}
		We can copy word by word the proof of \cite[Lemma 2.2.1]{caro2012overholonomicity}.
	\end{proof}
	
	\begin{lemma} \label{lem:log-D-push}
		In the situation of Lemma \ref{lem:log-D-dagger}, let $\wt{u} : \wt{\Xf}^\sharp \to \Xf$ and $u : \Xf^\sharp \to \Xf$ be the morphisms. Then for $\wt{\E} \in D^b_{\mathrm{coh}}(\D^\dag_{\wt{\Xf}^\sharp}(^\dag D)_\Qb)$ and $\E \in D^b_{\mathrm{coh}}(\D^\dag_{\Xf^\sharp}(^\dag D)_\Qb)$, we have the canonical isomorphisms:
		\begin{align*}
			u_{D \cup Z', +}(\wt{\E}(^\dag Z')) \simeq \wt{u}_{D \cup Z', +}(\wt{\E}(^\dag Z')) \simeq (\wt{u}_{D, +}(\wt{\E}))(^\dag Z'), \\
			\wt{u}_{D \cup Z', +}(\E(^\dag Z')) \simeq u_{D \cup Z', +}(\E(^\dag Z')) \simeq (u_{D, +}(\E))(^\dag Z').
		\end{align*}
	\end{lemma}
	\begin{proof}
		We can copy word by word the proof of \cite[Lemma 2.2.2]{caro2012overholonomicity}.
	\end{proof}
	
	\begin{proposition} \label{prop:rho}
		In the situation of {\rm \ref{p-smooth-retraction}}, let $D$ be a divisor of $X$, $\Xf^\sharp = (\Xf, \Zf)$ and $u : \Xf^\sharp \to \Xf$ be the canonical morphism. For a coherent $\D^\dag_{\Xf^\sharp}(^\dag D)_\Qb$-module locally projective of finite rank over $\O_\Xf(^\dag D)_\Qb$ satisfying the strong exponent condition, the canonical morphism $\rho : u_{D, +}(\E) \to \E(^\dag Z)$ is an isomorphism.
	\end{proposition}
	\begin{proof}
		The proof is similar to \cite[Proposition 2.2.3]{caro2012overholonomicity}. By Mayer-Vietoris exact triangle, it suffices to prove that $\RG^\dag_{Z_i} u_{D, +} \E = 0$ for each $i = 1, \dots, m.$ By \ref{p-smooth-retraction}, we have the commutative diagram
		\begin{align*}
			\xymatrix{
				\Zf_{i} \ar@{_(->}[rd]^-{\iota'} \ar@{_(->}@/_10pt/[rdd]_-{\iota} \ar@{=}@/^15pt/[rrd] &  &\\
				& \Xf' \ar[d]^-{f} \ar[r]^-{g} & \Zf_{i} \\
				& \Xf.
			}
		\end{align*}
		with $f$ affine etale, $f^{-1}(\Zf_i) = \iota'(\Zf_i)$ and $g$ with a finite $p$-basis.
		Let $\Zf' = f^{-1}(\Zf)$, $\Xf'^\sharp = (\Xf', \Zf')$, $D' = f^{-1}(D)$, $\Zf_i^\sharp = (\Zf_i,  \cup_{j \neq i} (\Zf_i \cap \Zf_j))$ and let $u' : \Xf'^\sharp \to \Xf'$, $v : \Zf_i^\sharp \to \Zf_i$, $f^\sharp : \Xf'^\sharp \to \Xf^\sharp$ and $g^\sharp : \Xf'^\sharp \to \Zf_i^\sharp$ be canonical morphisms. Then by Lemma \ref{lem:etale-pullback}, $f^{\sharp, !}_D\E$ is a coherent $\D^\dag_{\Xf'^\sharp}(^\dag D')_\Qb$-module locally projective of finite rank over $\O_{\Xf'}(^\dag D')_\Qb$. Using Corollary \ref{cor:comp-arithB}, $g_+ \RG^\dag_{Z_i} u'_{D', +} f^{\sharp, !}_D\E = 0.$ Using Lemma \ref{lem:u+comm}, we have $0 = g_+ \RG^\dag_{Z_i} u'_{D', +} f^{\sharp, !}_D\E = g_+ \iota'_+ \iota'^! u'_{D', +} f^{\sharp, !}_D \E = \iota'^! u'_{D', +} f^{\sharp, !}_D \E = \iota'^! f^{!}_D u_{D, +} \E = \iota^! u_{D, +} \E$. Hence $\RG^\dag_{Z_i} u_{D, +}\E = \iota_+ \iota^! u_{D, +} \E = 0.$
	\end{proof}
	
	\subsection{$p$-smooth devissage in isocrystals}
	\begin{definition} \label{def:p-sm-dev}
		Let $\Xf$ be a formal scheme over $\V[[t]]$ locally with finite $p$-bases over $\V$, $Y \subseteq X$, $T \subseteq X$ a divisor such that $Y = \ol{Y} - T$ and $\ol{Y}$ locally has finite $p$-bases over $k$. We say that a sequence of divisors $T_1, \dots, T_r \supseteq T$ of $X$ with $T_r = T$ gives a $p$-smooth devissage of $\E \in \LD{qc}(\wh{D}^{(\bullet)}_\Xf)$ over $Y$ in isocrystals with respect to $T$ if
		\begin{itemize}
			\item for $Y_i := \ol{Y} \cap T_1 \cap \cdots \cap T_i - T_{i+1}$, $\ol{Y_i}$ locally has finite $p$-bases over $k$ for each $i = 0, \dots, r-1$;
			\item For any $j$, $H^j(\RG^\dag_{Y_i}(\E))$ belongs to $\MIC^\ddag(\ol{Y_i}, \Xf, T_{i+1})$.
		\end{itemize}
		Note that $\RG^\dag_{Y_i} \E = \RG^\dag_{\ol{Y_i}} (^\dag T_{i+1}) \E.$ When there is some sequence of divisors giving a $p$-smooth devissage of $\E$ over $Y$ with respect to $T$, we simply say that $\E$ is $p$-smoothly devissable over $Y$ in isocrystals with respect to $T$. When $Y$ or $T$ is obvious from the context, we omit them.
	\end{definition}
	
	\begin{lemma}
		Let $\Xf$ be a formal scheme over $\V[[t]]$ locally with finite $p$-bases over $\V$, $Y \subseteq X$, $T \subseteq X$ a divisor such that $Y = \ol{Y} - T$ and $\ol{Y}$ locally has finite $p$-bases over $k$. Let $T_1, \dots, T_r \supseteq T$ be a sequence of divisors of $X$ with $T_r = T$ and $\E \in \LD{qc}(\wh{D}^{(\bullet)}_\Xf)$. Put $Y_i := \ol{Y} \cap T_1 \cap \cdots \cap T_i - T_{i+1}$ for $i = 0, \dots, r-1$. Suppose that each $\ol{Y_{i}}$ locally has finite $p$-bases over $k$ and that $\E$ is $p$-smoothly devissable in isocrystals over $Y_i$ with respect to $T_{i+1}.$ Then $\E$ is $p$-smoothly devissable in isocrystals over $Y$ with respect to $T.$
	\end{lemma}
	\begin{proof}
		We can copy word by word \cite[2.3.3]{caro2012overholonomicity}.
	\end{proof}
	
	\begin{lemma} \label{lem:ret-lem}
		Let $(\Xf, \Zf)$ be a strictly semistable pair over $\V[[t]]$ or a smooth pair over $\V$, $\Zf_1, \dots, \Zf_m$ the irreducible components of $\Zf$ and $D$ a divisor of $X$. Suppose that $Z_i$ is not contained in $D$ and that we have the commutative diagram as in {\rm (\ref{diag:ret})}
		\begin{align*}
			\xymatrix{
				\Zf_{i} \ar@{_(->}[rd]^-{\iota'} \ar@{_(->}@/_10pt/[rdd]_-{\iota} \ar@{=}@/^15pt/[rrd] &  &\\
				& \Xf' \ar[d]^-{f} \ar[r]^-{g} & \Zf_{i} \\
				& \Xf.
			}
		\end{align*}
		with $f$ affine etale, $f^{-1}(\Zf_i) = \iota'(\Zf_i)$ and $g$ with a finite $p$-basis. Let $\Zf^\sharp = (\Zf_i, \cup_{j \neq i} \Zf_i \cap \Zf_j)$, $\Xf'^\sharp = (\Xf', f^{-1}(\Zf))$ and let $f^\sharp : \Xf'^\sharp \to \Xf^\sharp$, $g : \Xf'^\sharp \to \Zf_i^\sharp$, $v : \Zf_i^\sharp \to \Zf_i$, $u' : \Xf'^\sharp \to \Xf'$ be the canonical morphism. Then for a bounded complex $\E$ of coherent $\D^\dag_{\Xf^\sharp}(^\dag D)$-modules locally projective of finite rank over $\O_{\Xf}(^\dag D)_\Qb$ satisfying exponent condition, $\iota^! u_{D, +}\E \simeq v_{D \cap Z_i, +} \G$ for some bounded complex $\G$ of coherent $\D^\dag_{\Zf_i^\sharp}(^\dag D\cap Z_i)$-modules locally projective of finite rank over $\O_{\Zf_i}(^\dag D \cap Z_i)_\Qb$ satisfying exponent condition.
	\end{lemma}
	\begin{proof}
		Let $D' = f^{-1}(D).$ By Lemma \ref{lem:u+comm}, Lemma \ref{lem:u+RGcom} we have
		$\iota^! u_{D, +} \E = \iota'^! f^! u_{D, +} \E \overset{\ref{lem:u+comm}}{=} \iota'^! u'_{D', +} f^{\sharp, !} \E = g_+ \iota'_+ \iota'^! u'_{D', +} f^{\sharp, !} \E = g_+ \RG^\dag_{Z_i} u'_{D', +} f^{\sharp, !}\E \overset{\ref{lem:u+RGcom}}{=} g_+u'_{D', +} \RG^\dag_{Z_i}f^{\sharp, !}\E = v_{D \cap Z_i, +}g^\sharp_+\RG^\dag_{Z_i}f^{\sharp, !}\E.$ Noting that $f^{\sharp, !}$ is exact, we can take $\G = g^\sharp_+\RG^\dag_{Z_i}f^{\sharp, !}\E$ by Lemma \ref{lem:etale-pullback} and Theorem \ref{thm:comp-arithB}. 
	\end{proof}
	
	\begin{theorem}
		$(\Xf, \wt{\Zf})$ be an affine strictly semistable pair over $\V[[t]]$ with a finite $p$-basis over $\V$ or a smooth pair over $\V$ with an etale coordinate over $\V$ and let $\Zf_0, \dots, \Zf_m$ be the irreducible components of $\wt{\Zf}$. Let $I, I'$ be disjoint subsets of $\{ 0, \dots, m \}$, $\Zf := \cup_{i \in I} \Zf_i$, $\Zf' := \cup_{i \in I'} \Zf'_i$ and let $\Xf^\sharp = (\Xf, \tilde{\Zf})$. Then there exists a sequence of subdivisors $T_0, \dots, T_N$ of $Z$ with the following property: each $T_i$ is a subdivisor of $\Zf$ and for any bounded complex $\E$ consisting of coherent $\D^\dag_{\Xf^\sharp}(^\dag Z')_\Qb$-modules locally projective of finite rank over $\O_{\Xf}(^\dag Z')_\Qb$ satisfying the exponent condition, $T_0 \cup Z', \dots, T_N \cup Z'$ gives a $p$-smooth devissage of $\F := u_{Z', +}\E$ over $X - Z'$ in isocrystals.
	\end{theorem}
	\begin{proof}
		When $I$ is empty, there is nothing to prove; by reordering we may assume $I = \{ 0, \dots, r \}$ and $I' = \{ r+1, \dots, m'\}$. For a binary expansion $\displaystyle k = \sum_{i=0}^r \epsilon_i 2^i$ of $0 \leq k \leq 2^{r+1}-1,$ put $\displaystyle T_k = \bigcup_{\epsilon_i = 0} Z_i$ and $T'_k = T_k \cup Z'.$ We will prove that $T'_0, \dots, T'_{2^{r+1}-1}$ gives $p$-smooth devissage of $\F$ over $X - Z'$ in isocrystals by induction on $r$. Let $Y_k = T'_0 \cap \cdots \cap T'_{k-1} - T'_k$. We have to prove that $\ol{Y_k}$ locally has finite $p$-bases over $k$ and that $H^j(\RG^\dag_{Y_k} \F)$ belongs to $\MIC^\ddag(\ol{Y_k}, \Xf, T'_{k} / K)$. 
		
		(i) $0 \leq k \leq 2^r - 1$: By induction, the statement applied with $I = \{ 0, \dots, r-1 \}$ and $I' = \{ r, \dots, m' \}$ to $\E(^\dag Z_r)$ asserts that for the binary expansion $k = \displaystyle\sum_{i = 0}^{r-1} \epsilon_i 2^i$, put $\displaystyle S_k = \bigcup_{\epsilon_i = 0} T_i$, $S'_k = S_k \cup Z_r \cup Z'$ and $Y'_k = S'_0 \cap \cdots \cap S'_{k-1} - S'_k$; then  $\ol{Y'_k}$ locally has finite $p$-bases over $k$ and $H^j(\RG^\dag_{Y_k} u_{Z_r \cup Z', +}(\E(^\dag Z_r)))$ belongs to $\MIC^\ddag(\ol{Y_k}, \Xf, S'_k / K).$ Since $T_k = S_k \cup Z_r$ and $T'_k = S'_k$, we have $Y_k = Y'_k$ and hence $\ol{Y_k} = \ol{Y'_k}$ locally has finite $p$-basis over $k$. Moreover, $\RG^\dag_{Y'_k} u_{Z_r \cup Z', +}(\E(^\dag Z_r)) = \RG^\dag_{Y_k}(u_{Z', +}(\E)(^\dag Z_r)) = \RG^\dag_{Y_k}(\F)$ by Lemma \ref{lem:log-D-push}. Note that $Z_r$ and $Y_k$ are disjoint and therefore $\RG^\dag_{Y_k}(^\dag Z_r) = \RG^\dag_{Y_k}$.
		
		(ii) $2^r \leq k \leq 2^{r+1} - 1$: By \ref{p-smooth-retraction}, we have the commutative diagram
		\begin{align*}
			\xymatrix{
				\Zf_{r} \ar@{_(->}[rd]^-{\iota'} \ar@{_(->}@/_10pt/[rdd]_-{\iota} \ar@{=}@/^15pt/[rrd] &  &\\
				& \Xf' \ar[d]^-{f} \ar[r]^-{g} & \Zf_{r} \\
				& \Xf.
			}
		\end{align*}
		with $f$ affine etale, $f^{-1}(\Zf_r) = \iota'(\Zf_i)$ and $g$ with a finite $p$-basis. Let $v : \Zf^\sharp_r = (\Zf_r, \cup_{j \neq r} \Zf_r \cap \Zf_j) \to \Zf$ be the canonical morphism. By Lemma \ref{lem:ret-lem}, we have a bounded complex $\G$ of coherent $\D^\dag_{\Zf^\sharp_r}(^\dag Z_r \cap Z')_\Qb$-modules locally projective of finite rank over $\O_{\Zf_r}(^\dag Z_r \cap Z')_\Qb$ satisfying the exponent condition such that $\iota^! u_{Z', +} \E \simeq v_{Z_r \cap Z', +}\G$. By inductive hypothesis, for $0 \leq k' \leq 2^r -1$ and its binary expansion $\displaystyle \sum_{i = 0}^{r-1} \epsilon'_i 2^i$, put $S_{k'} = \displaystyle\bigcup_{\epsilon_i = 0} (Z_i \cap Z_r)$, $S'_{k'} = S_{k'} \cup (Z_r \cap Z')$ and $Y'_{k'} = S'_0 \cap \cdots \cap S'_{k'-1} - S'_{k'} \subseteq Z_r$; then $\ol{Y'_{k'}}$ locally has finite $p$-bases over $k$ and $H^j(\RG^\dag_{Y'_{k'}} v_{Z_r \cap Z', +}\G)$ belongs to $\MIC^\ddag(\ol{Y'_{k'}}, \Zf_r, S'_{k'} / K).$ Since $T'_0 \cap \cdots \cap T'_{2^r} = Z_r \cup Z'$ and $Z'_r \cap T'_k = S'_{k-2^r}$, we have $Y_k = (Z_r \cup Z') \cap T'_{2^r+1} \cap \cdots \cap T'_{k-1} - T'_k = S'_0 \cap \cdots \cap S'_{k-2^r-1} - S'_{k-2^r} = Y'_{k-2^r}.$ Therefore $\ol{Y_k} = \ol{Y'_{k-2^r}}$ locally has finite $p$-basis. Moreover $\iota_+ \RG^\dag_{Y'_{k-2^r}} v_{Z_r \cap Z', +} \G = \iota_+ \RG^\dag_{Y'_{k-2^r}} \iota^! u_{Z', +} \E = \iota_+ \iota^! \RG^\dag_{Y_k} u_{Z', +} \E = \RG^\dag_{Y_k} u_{Z', +} \E$ as $Y_k \subseteq Z_r$. By Proposition \ref{prop:sp+-pullback}, noting that $\iota_+$ is a quasi-inverse to $\iota^*$, $H^j(\RG^\dag_{Y_k} u_{Z', +} \E)$ belongs to $\MIC^\ddag(\ol{Y_k}, \Xf, T'_k / K).$
	\end{proof}
	
	\begin{remark}
		We assumed that $(\Xf, \Zf)$ has a global finite $p$-basis in order to invoke \ref{p-smooth-retraction}. However, the sequence of divisors constructed in the proof using binary expansions keeps the validity of statement without global finite $p$-basis, since being an overconvergent isocrystal (more precisely, in the essential image of $\sp_+$ functors) is a local property on $\Xf$.
	\end{remark}
	
	\subsection{Overholonomicity of log extendable isocrystal}
	\begin{proposition} \label{prop:log-relative-duality}
		Let $(\Xf, \Zf)$ be an affine strictly semistable pair over $\V[[t]]$, $u : \Xf^\sharp \to \Xf$ the canonical morphism and $\E$ a coherent $\D^\dag_{\Xf^\sharp, \Qb}$-module locally projective of finite rank over $\O_{\Xf, \Qb}$. Then we have a canonical isomorphism $\Db u_+ \E \simeq u_+ (\E^\vee(-\Zf))$, where $\E^\vee$ is the dual of corresponding log isocrystal.
	\end{proposition}
	\begin{proof}
		The proof in the classical case works as well, see \cite[2.3.12]{caro2012overholonomicity}.
	\end{proof}
	
	\begin{theorem} \label{thm:u+-ovhol}
		Let $(\Xf, \Zf)$ be a strictly semistable pair over $\V[[t]]$ or a smooth pair over $\V$, $\E$ a complex of coherent $\D^\dag_{\Xf^\sharp, \Qb}$-modules locally projective of finite rank over $\O_{\Xf, \Qb}$ satisfying exponent condition. For $u : \Xf^\sharp \to \Xf$ the canonical morphism, $u_+\E$ is overholonomic.
	\end{theorem}
	\begin{proof}
		Consider following statements:
		\begin{itemize}
			\item $(P_{n, r})$: given $(\Xf, \Zf)$ and $\E$ as above, $u_+ \E$ is $(n, r)$-overholonomic.
			\item $(Q_{n, r})$: given $(\Xf, \Zf)$ and $\E$ as above, $\RG^\dag_Z u_+ \E$ is $(n, r)$-overholonomic.
			\item $(R_{n, r})$: given $(\Xf, \Zf)$ and $\E$ as above, $(^\dag Z)u_+ \E = (^\dag Z) \E$ is $(n, r)$-overholonomic.
		\end{itemize}
		
		(i) First we prove that for $n \geq 1$ and $r \geq -1$, $(P_{n-1,r})$ implies $(Q_{n, r})$.
		Let $\Zf_1, \dots, \Zf_m$ be the irreducible components of $\Zf$. By Mayer-Vietoris exact triangle, it suffices to prove that each $\RG^\dag_{Z_i} u_+ \E$ is $(n, r)$-overholonomic. The assertion is local on $\Xf$, so we may assume that we are in the situation of \ref{p-smooth-retraction}. We have the commutative diagram
		\begin{align*}
			\xymatrix{
				\Zf_{i} \ar@{_(->}[rd]^-{\iota'} \ar@{_(->}@/_10pt/[rdd]_-{\iota} \ar@{=}@/^15pt/[rrd] &  &\\
				& \Xf' \ar[d]^-{f} \ar[r]^-{g} & \Zf_{i} \\
				& \Xf.
			}
		\end{align*}
		with $f$ affine etale, $f^{-1}(\Zf_i) = \iota'(\Zf_i)$ and $g$ with a finite $p$-basis. Let $v : \Zf^\sharp_i = (\Zf_i, \cup_{j \neq i} \Zf_r \cap \Zf_j) \to \Zf$ be the canonical morphism. By Lemma \ref{lem:ret-lem}, we have a bounded complex $\G$ of coherent $\D^\dag_{\Zf^\sharp_r, \Qb}$-modules locally projective of finite rank over $\O_{\Zf_r, \Qb}$ satisfying the exponent condition such that $\iota^! u_{+} \E \simeq v_{+} \G$. By inductive hypothesis, $v_+\G$ is $(n-1, r)$-overholonomic; by Proposition \ref{prop:nr-ovhol}-(i), $\iota_+v_{+} \G \simeq \RG_{Z_i}^\dag u_{+} \E$ is $(n, r)$-overholonomic.
		
		(ii) Next we show that for $n \geq 0$ and $r \geq 0$, $(P_{n, r-1})$ and $(Q_{n, r})$ implies $(R_{n, r})$. By devissage, we may assume that $\E$ is a module. When $r = 0$ (resp. $r \geq 1$), we have to check that $(^\dag Z)\E$ is $(n, r-1)$-overholonomic and that for any smooth morphism $f : \Xf' \to \Xf$ with $\delta_{X'} \leq n$ and any divisor $T'$ of $X'$, the complex $(^\dag T')f^! (^\dag Z)\E$ (resp. $\Db(^\dag T')f^! (^\dag Z)\E$) is $(n, r-1)$-overholonomic. Let $f^\sharp : \Xf'^\sharp = (\Xf', f^{-1}(Z)) \to \Xf$ be the canonical morphism. Noting that $f^{\sharp, !}[\delta_X - \delta_{X'}]$ is exact, by Lemma \ref{lem:etale-pullback}, $f^{\sharp, !} \E$ is a complex over $\Xf'$ satisfying the hypothesis of this theorem, and $(^\dag T')f^!(^\dag Z) \E = (^\dag T')(^\dag f^{-1}(Z)) f^{\sharp, !} \E$. By replacing $\Xf'^\sharp$ by $\Xf^\sharp$ and $\E$ by $f^{\sharp, !}\E$, we are reduced to prove that for any divisor $D$ of $X$, $(^\dag D \cup Z)\E$ is (resp. $(^\dag Z)\E$ and $\Db(^\dag D \cup Z)\E$ are) $(n, r-1)$-overholonomic.
		
		Henceforth we suppose that $(\Xf, \Zf)$ is a strictly semistable pair over $\V[[t]]$; the case of smooth pairs can be proved in the similar way, or just refer to \cite[Theorem 2.3.13]{caro2012overholonomicity}. By de Jong's alteration theorem, we have a projective alteration $a_0 : Y \to X$ over $k[[t]]$ and a finite extension $l[[u]]$ of $k[[t]]$ such that $(Y, a_0^{-1}(D \cup Z)_{\red})$ is a strictly semistable pair over $l[[u]]$. There is $n \geq 0$ such that $l[[u]]$ is etale over $k[[t^{1/p^n}]]$. Let $\W[[u]]$ be the finite etale extension of $\V[[t^{1/p^n}]]$ lifting $l[[u]]$ over $k[[t^{1/p^n}]]$. Let $p : \Xf' := \Xf \times_{\Spf \V[[t]]} \Spf \W[[l]] \to \Xf^{(n)} := \Xf \times_{\Spf \V[[t]]} \Spf \V[[t^{1/p^n}]]$ and let $q : \Xf^{(n)} \to \Xf$ be canonical projections and let $a'_0 : Y \to X'$ be the canonical map.
		We can choose some formal projective space $f : \Pf = \wh{\mathbb{P}}^d_{\Xf'} \to \Xf'$ and a closed immersion $Y \subseteq P'$ whose composition with $f : P' \to X'$ is $a'_0 : Y \to X'$. Since $\E' := (qp)^! \E (^\dag D \cup Z)$ corresponds to an overconvergent isocrystal over $(X' - p^{-1}(D \cup Z), X', \Xf')/K$, $\E'$ is a direct summand of $f_+\RG^\dag_{Y} f^! \E'$ by Proposition \ref{prop:dir-sum-isoc}. Since $p$ is finite etale, $q_+q^!\E(^\dag D \cup Z)$ is a direct summand of $q_+p_+p^!q^! \E(^\dag D \cup Z) = (qp)_+ \E'$ and since $q_0 : X^{(n)} \to X$ is finite surjective and radicial, by \cite[Proposition 11.4.6]{caro2021arithmetic}, $\E(^\dag D \cup Z) \simeq q^!q_+ \E(^\dag D \cup Z)$; thus $\E$ is a direct summand of $(qp)_+ \E'$. Again by Proposition \ref{prop:nr-ovhol}-(i), it suffices to prove that $\RG^\dag_Y f^! \E'$ is (resp. $\RG^\dag_Y f^! \E'$ and $\Db\RG^\dag_Y f^! \E'$ are) $(n+d, r-1)$-overholonomic.
		
		The assertion is local $\Pf$ and we may assume that $\Pf$ is affine and that there is a lift $\iota : \Yf \hookrightarrow \Pf$ of the closed immersion over $\W$ and a lift $\Zf'$ of $Z' := a_0'^{-1}(D \cup Z)_{\red}$ such that $(\Yf, \Zf')$ is a strictly semistable pair over $\W[[u]]$. More precisely, when there is an etale morphism $Y \to \Spec l[[u]][t_1, \dots, t_n]/(u - t_1 \cdots t_r)$ and $Z'$ is cut out by $t_1 \cdots t_m$, choose an etale lift $\Yf \to \Spf \W[[u]] \allowbreak \langle t_1, \dots, t_n \rangle / (u - t_1 \cdots t_r)$ (which is unique by henselian property). Since $\Pf$ is formally etale over $\V$, one gets a lift $\Yf \hookrightarrow \Pf$ over $\W$. Put $a = f \iota : \Yf \to \Xf'$ and $u : \Yf^\sharp \to \Yf$. Since $\RG^\dag_{X'}f^! \E' = \iota_+ \iota^! f^! \E' = \iota_+ a^! \E'$, again by Proposition \ref{prop:nr-ovhol}-(i), it suffices to prove that $a^! \E'$ is (resp. $a^! \E'$ and $\Db a^! \E'$ are) $(n, r-1)$-overholonomic.
		
		Since $a_0$ is flat, so is $a$ flat; therefore $a^! \E'$ is, up to shift, isomorphic to the pullback $a^*\E'$ of overconvergent isocrystals. Moreover, for $(qpa)^\sharp : \Yf^\sharp \to \Xf^\sharp$ the canonical morphism, we have $a^* \E' \simeq (^\dag Z')(qpa)^{\sharp, *} \E$. Note that here $(qpa)^{\sharp, *}$ is the functor from $\Isoc(\Xf^\sharp/\V)$ to $\Isoc(\Yf^\sharp/\V)$ and $(\Yf, \Zf')$ is not necessarily strictly semistable over $\V[[t]]$ (it is over $\W[[u]]$); however as $\W$ is etale over $\V$, we have canonical equivalence $\Isoc(\Yf^\sharp/\V) \simeq \Isoc(\Yf^\sharp/\W)$ and an isomorphism $\D^\dag_{\Yf^\sharp/\V} \simeq  \D^\dag_{\Yf^\sharp/\W}$, so one still have the equivalence of Proposition \ref{prop:isoc-Dmod-equiv}. Therefore we have to prove that for $\F = a^{\sharp, *} (qp)^{\sharp, *} \E$, $(^\dag Z') \F$ is (resp. $(^\dag Z')\F$ and $\Db (^\dag Z')\F$ are) $(n, r-1)$-overholonomic. Note that $\F$ is a module over $\Yf^\sharp$ satisfying the hypothesis of this theorem; by exactness of $\Db$ and $(^\dag Z')$, by devissage we may assume that $\F$ is a module concentrated at degree zero.
		
		When $r = 0$, since $\F$ over $\Yf^\sharp$ satisfies the hypothesis of this theorem, $(Q_{n, 0})$ implies that $\RG^\dag_{Z'}u'_+ \F$ is $(n, 0)$-overholonomic; hence $(n, -1)$-overholonomic. $(P_{n, -1})$ implies that $u'_+ \F$ is $(n, -1)$-overholonomic. By devissage, $(^\dag Z') \F = (^\dag Z') u'_+ \F$ is $(n, -1)$-overholonomic.
		
		When $r \geq 1$, by similar argument, $(^\dag Z') \F$ is $(n, r-1)$-overholonomic. Since $\RG^\dag_{Z'}u'_+ \F$ is $(n, r)$-overholonomic by $(Q_{n, r})$, $\Db \RG^\dag_{Z'}u'_+ \F$ is $(n, r-1)$-overholonomic. $\Db u'_+ \F \simeq u'_+ (\F(-\Zf'))$ by Proposition \ref{prop:log-relative-duality} and it is $(n, r-1)$-overholonomic by $(P_{n, r-1})$. By the exact triangle
		$$ \Db (^\dag Z') \F \to \Db u'_+\F \to \Db \RG^\dag_{Z'} u'_+ \F \to,$$
		we conclude that $ \Db (^\dag Z') \F$ is $(n, r-1)$-overholonomic.
		
		(iii) We have $(P_{n-1, r}) \wedge (P_{n, r-1}) \Rightarrow (Q_{n, r}) \wedge (P_{n, r-1}) \Rightarrow (Q_{n, r}) \wedge (R_{n, r}) \Rightarrow (P_{n, r})$. $(P_{0, r})$ is obviously true, since if $\delta_{\Xf} = 0$, $\Xf = \Spf \W$ is finite etale over $\Spf \V$ and $D^b_{\rm coh}(\D^\dag_{\Xf, \Qb}) = D^b_{\rm fin}(L)$ with $L$ the quotient field of $\W$. $(P_{n, -1})$ comes from (\ref{eq:u+}).
	\end{proof}
	
	\section{Six functor formalism of arithmetic $\D$-modules}
	\subsection{Frobenius structures}
	\begin{_empty}
		$\Pf$ be a formal $\V[[t]]$-scheme of finite type locally having finite $p$-bases over $\V$. Recall that we have a lift $\sigma : \V \simto \V$ of the $s$-th Frobenius $k \to k, \ x \to x^{p^s}$. Then by \cite[Corollary 4.5.3]{caro2021arithmetic}, we have the exact $\sigma$-linear Frobenius pullback functors
		$$F^* : D_{\rm (coh)}(\D^\dag_{\Pf, \Qb}) \simto D_{\rm (coh)}(\D^\dag_{\Pf, \Qb})$$
		and
		$$F^* :\LD{?}(\wh{\D}^{(\bullet)}_{\Pf}) \simto \LD{?}(\wh{\D}^{(\bullet)}_{\Pf})$$
		for $?$ either ${\rm qc}$ or ${\rm coh}$. The Frobenius pullback functor preserve various categories $\MIC^{(\bullet)}$ and $\MIC^\ddag$.
		$F^*$ is compatible with various operations such as extraordinary pullbacks, pushforwards and others.
		An $a$-th Frobenius structure of $\E \in D(\D^\dag_{\Pf, \Qb})$ (resp. $\E \in \LD{qc}(\wh{\D}^{(\bullet)}_{\Pf})$) is an isomorphism $\phi : \E \to (F^*)^a \E$.
		
		Let $f : \V \to \W$ be a morphism of the form $\V \to \V \otimes_{W(k)} W(l)$ for $l$ a perfect field over $k$. Then $\sigma$ extends to a lift $\sigma' : \W \to \W$ of the $s$-th Frobenius on $l$ and we can discuss the Frobenius structures over $\W$ as well. Furthermore, two Frobenius pullbacks are compatible under the base change functor $f^*$.
		
		Let $X$ be a closed subscheme of $P$ locally having finite $p$-basis over $k$ and $T$ a divisor of $P$ such that $T \cap X$ is a divisor of $X$. Similarly, we have the $\sigma$-linear Frobenius pullback functor
		$$F^* : \MIC^\dag(X, \Pf, T / K) \to \MIC^\dag(X, \Pf, T / K).$$
		One easily checks that two functors $F^*$ are compatible through the functor $\sp_{X \hookrightarrow \Pf, T, +}$. Thus, by construction they commute with general $\sp_{X \hookrightarrow \Pf, Z, +}$ for any closed subschemes $X, Z$ not necessarily having finite $p$-bases over $k$.
	\end{_empty}
	
	\begin{theorem} \label{thm:isoc-ovhol}
		Let $\Pf$ be a formal $\V[[t]]$-scheme of finite type locally having finite $p$-bases over $\V$, $X, Z$ closed subschemes of $P$ and $Y := X - Z$ the subscheme locally having finite $p$-bases over $k$. Then for any $E \in \MIC^\dag(Y, X, \Pf / K)$ which can be endowed with some $a$-th Frobenius structure, the functor $\sp_{X \hookrightarrow \Pf, Z, +} : \MIC^\dag(Y, X, \Pf / K) \to \MIC^{(\bullet)}(X, \Pf, Z / K)$ maps $E$ into $\LD{oh}(\wh{\D}^{(\bullet)}_\Pf)$.
	\end{theorem}
	\begin{proof}
		Since Frobenius pullback is compatible with base changes, it suffices to prove that it maps $E$ into $\LD{ovhol}(\wh{\D}^{(\bullet)}_\Pf).$ As in \cite[Theorem 18.3.2.2]{caro2025arith}, we may assume that $\Pf$ is affine and connected and $Z = T$ is a divisor of $P$. By \cite[Theorem 7.4.7]{liu2026semistable}, there is a finite extension $k'[[t']]$ of $k[[t]]$, a projective alteration $a_0 : X' \to X$ such that for $Z' := a_0^{-1}(T)_\red$, $(X', Z')$ is a strictly semistable pair over $k'[[t']]$, and the isocrystal $E' := a_0^*E$ over $(X'-Z', X')/K'$ extends to a log convergent isocrystal with nilpotent residues over $Y^\sharp / K'$. Here we let $\V'[[t']]$ be the etale extension of $\V[[t^{1/p^n}]]$ corresponding to the etale extension $k'[[t']]$ of $k[[t^{1/p^n}]]$ and let $K'$ be the quotient field of $\V'$. Choose a closed immersion $X' \hookrightarrow \Pf'$ into a formal $\V'[[t']]$-scheme and a projective morphism $f : \Pf' \to \Pf$ over $\V[[t]]$ such that
		\[\xymatrix{
			X' \ar[r] \ar[d]_{a_0} & \Pf' \ar[d]^{f} \\
			X \ar[r] & \Pf
		}\]
		is commutative. For $Y' = X' - Z'$, $f$ induces the functor $f^*_K : \MIC^\dag(Y, X, \Pf / K) \to \MIC^\dag(Y', X', \Pf' / K) \allowbreak \simeq \MIC^\dag(Y', X', \Pf' / K')$ and $f^+ : \MIC^{(\bullet)}(X, \Pf, Z / K) \to \MIC^{(\bullet)}(X', \Pf', Z' / K) \allowbreak \simeq \MIC^{(\bullet)}(X', Z', \Pf' / K')$. Let $\E = \sp_{X \hookrightarrow \Pf, Z} (E)$; then $\RG^\dag_{Y'} f^! \E \simeq \sp_{X' \hookrightarrow \Pf', Z'} (a_0^* E)$. As in \cite[Lemma 16.1.11.2]{caro2025arith}, $\E$ is a direct summand of $f_+ \RG^\dag_{Y'} f^! \E$. Thus it suffices to prove that $\sp_{X' \hookrightarrow \Pf', Z'} (a_0^* E)$ is overholonomic. The assertion is local on $\Pf$, so we may assume that $\Pf$ is affine and that there exists a lift $\Xf'$ of $X$ to a formal $\V[[t]]$-scheme locally having finite $p$-bases over $\V$, a closed subscheme $\Zf'$ of $\Xf'$ so that $(\Xf', \Zf')$ is a strictly semistable pair over $\V$ and a lift $\iota : \Xf' \hookrightarrow \Pf'$ over $\V$ of the closed immersion $X' \hookrightarrow \Pf'$. Then $\sp_{X' \hookrightarrow \Pf', Z'} (a_0^* E) = \iota_+ \sp_{X' \hookrightarrow \Xf', Z'} (a_0^* E)$, and it suffices to prove that $\E' := \sp_{X' \hookrightarrow \Xf', Z'} (a_0^* E)$ is overholonomic. Note that $X'$ locally has finite $p$-bases over $k$, we can consider $\MIC^\ddag(X', \Xf', Z' / K)$ instead of $\MIC^{(\bullet)}(X', \Xf', Z' / K).$ Since $a_0^* E$ extends to a log convergent isocrystal with nilpotent residues, there is a coherent $\D^\dag_{\Xf'^\sharp, \Qb}$-module $\F$ locally projective of finite rank over $\O_{\Xf', \Qb}$ satisfying strong exponents conditions and $(^\dag Z')\F = \E'$. By Proposition \ref{prop:rho} and Theorem \ref{thm:u+-ovhol}, $\E' = u_+ \F$ is overholonomic.
	\end{proof}
	
	\begin{theorem} \label{thm:ovhol-dev}
		Let $\Pf$ be a formal $\V[[t]]$-scheme of finite type locally having finite $p$-bases over $\V$, $X, Z$ closed subschemes of $P$, $Y := X - Z$ the subscheme locally closed subscheme of $P$ locally having finite $p$-bases over $k$ and let $\E \in \LD{qc}(\wh{\D}^{(\bullet)}_\Pf)$ be a complex such that there is an isomorphism $\RG^\dag_Y \E \simeq \E$ and that can be endowed with some $a$-th Frobenius structure. Then the following conditions are equivalent:
		\begin{enumerate}[(i)]
			\item $\E \in \LD{povcoh}(X, \Pf, Z / \V)$;
			\item $\E \in \LD{ovcoh}(\wh{\D}^{(\bullet)}_\Pf)$;
			\item $\E \in \LD{ovhol}(\wh{\D}^{(\bullet)}_\Pf)$;
			\item $\E \in \LD{oh}(\wh{\D}^{(\bullet)}_\Pf)$;
			\item $\E \in \LD{dev}(Y, \Pf / \V)$.
		\end{enumerate}
	\end{theorem}
	\begin{proof}
		(iv) $\Rightarrow$ (iii) $\Rightarrow$ (ii) $\Rightarrow$ (i) are obvious. (i) $\Rightarrow$ (v) is Proposition \ref{prop:povcoh-dev}. Lastly, (v) implies (iv) by devissage and Theorem \ref{thm:isoc-ovhol}.
	\end{proof}

	\subsection{Construction of six functor formalism}
	\begin{_empty}
		Let $(X, Y, \Pf)$ be a formal frame with $\Pf$ of finite type over $\V[[t]]$ and locally having finite $p$-bases over $\V$. Let $D^b_{\oh}(X, Y, \Pf / K)$ be the triangulated full subcategory of $D^b_{\oh}(\D^\dag_{\Pf, \Qb})$ consisting of those $\E$ with some isomorphism $\RG^\dag_X \E \simeq \E.$
	\end{_empty}
	\begin{lemma} \label{lem:indep-frame}
		Let $(\id, f, u) : (X, Y', \Pf') \to (X, Y, \Pf)$ be a morphism of formal frames with $\Pf$ and $\Pf'$ quasi-projective smooth over $\V[[t]]$ and $f$ proper. Then the functors $\RG^\dag_X u^!$ and $u_+$ induces equivalences between $\D^b_{\oh}(X, Y, \Pf/K)$ and $\D^b_{\oh}(X, Y', \Pf'/K).$
	\end{lemma}
	\begin{proof}
		By \cite[Theorem 14.2.2]{caro2021arithmetic}.
	\end{proof}
	
	\begin{_empty} \label{emp:realizable}
		We call a scheme $X$ over $k[[t]]$ realizable if there is a formal frame $(X, Y, \Pf)$ with $\Pf$ quasi-projective and smooth over $\V[[t]]$ and $Y$ proper over $k[[t]]$. By Lemma \ref{lem:indep-frame}, it makes sense to define $D^b_{\oh}(X/K)$ for a realizable scheme $X$. One also has the Frobenius map $F^* : D^b_{\oh}(X/K) \simto D^b_{\oh}(X/K).$ We define the category of overholonomic arithmetic $\D$-modules of Frobenius type $D^b_{\oh, F}(X/K)$ to be the triangulated full subcategory of $D^b_{\oh}(X/K)$ generated by those objects $\E$ which can be endowed with some $a$-th Frobenius structure $F^{a*}\E \simto \E.$
		
		Similarly for a realizable scheme $X$ it makes sense to define $\MIC^\dag(X/K)$ as $\MIC^\dag(X, Y, \Pf/K)$ for some enclosing frame. By the equivalence $\Isoc(X, Y, \Pf/K) \simeq \MIC^\dag(X, Y, \Pf/K)$, we prefer to denote it by $\Isoc(X/K)$. Define $\Isoc_F(X/K)$ to be the thick full subcategory of $\Isoc(X/K)$ generated by those $E$ which can be endowed with some $a$-th Frobenius structure $F^{a*}E \simto E$. Choose an enclosing frame $(X, Y, \Pf)$. Theorem \ref{thm:isoc-ovhol} implies that for any closed subscheme $Z$ of $P$ such that $X = Y-Z$, we have the fully faithful functor $\sp_{Y \hookrightarrow \Pf, Z, +} : \Isoc_F(X, Y, \Pf / K) \to D^b_{\oh, F}(X, Y, \Pf/K).$ This is independent of the choice of frames or $Z$ by Lemma \ref{lem:indep-frame} and Proposition \ref{prop:bullet-indep-frame} and we get the functor $\sp_X : \Isoc_F(X/K) \to D^b_{\oh, F}(X/K).$
	\end{_empty}
	
	\begin{definition}[six operations]\label{def:six}
		Let $X$ be a realizable $k[[t]]$-scheme. An enclosing frame of $X$ is a formal frame $(X, Y, \Pf)$ with $\Pf$ projective and smooth over $\V[[t]]$ and $Y$ proper over $k[[t]].$ Let $f : X' \to X$ be a morphism of realizable $k[[t]]$-schemes. An enclosing morphism of $f$ is a morphism of enclosing frames $(f, g, u) : (X', Y', \Pf') \to (X, Y, \Pf)$ with $g$ proper. Note that it is always possible to take some morphism of enclosing frames.
		\begin{itemize}
			\item For a morphism $f : X' \to X$ of realizable schemes, choose an enclosing morphism $(f, g, u) : (X', Y', \Pf') \to (X, Y, \Pf)$. Define the functor
			$$f^! : D^b_{\oh, F}(X/K) \to D^b_{\oh, F}(X'/K)$$
			to be $\RG^\dag_{X'} u^! (-)$.
			\item For a morphism $f : X' \to X$ of realizable schemes, choose an enclosing morphism $(f, g, u) : (X', Y', \Pf') \to (X, Y, \Pf)$. Define the functor
			$$f^+ : D^b_{\oh, F}(X'/K) \to D^b_{\oh, F}(X/K)$$
			to be $u_+(-)$.
			\item For a realizable scheme $X$, choose a morphism of enclosing frame $(X, Y, \Pf)$. Define the functor 
			$$\Db_X : D^b_{\oh, F}(X/K)^\op \to D^b_{\oh, F}(X/K)$$
			to be $\RG^\dag_X \Db_{\Pf}(-).$
			\item For a morphism $f : X' \to X$ of realizable schemes, define the morphisms 
			$$f^+ : D^b_{\oh, F}(X/K) \to D^b_{\oh, F}(X'/K)$$
			and
			$$f^! : D^b_{\oh, F}(X'/K) \to D^b_{\oh, F}(X/K)$$
			by $f^+ = \Db_{X'} f^! \Db_X$ and $f_! = \Db_{X} f_+ \Db_{X'}.$
			\item For a realizable scheme $X$, choose a morphism of enclosing frame $(X, Y, \Pf)$. Define the functor
			$$(-) \wt{\otimes}_X (-) : D^b_{\oh, F}(X/K) \times D^b_{\oh, F}(X/K) \to D^b_{\oh, F}(X/K)$$
			to be $(-) \, \wt{\otimes}^{\mathbb{L}^\dag}_{\O_\Pf} (-) [-\delta_P]$. Also define $(-) \otimes_X (-) = \Db(\Db(-) \wt{\otimes}_X \Db(-)).$
		\end{itemize}
	\end{definition}
	
	\begin{remark}
		To check that these functors are well-defined, one has to check  that they are independent of the choice of enclosing frames or enclosing morphisms and that they send objects of $D^b_{\oh, F}$ into $D^b_{\oh, F}$. Independence can be checked as in \cite[\S 12.2]{caro2019arithmetic}. For preservation of $D^b_{\oh, F}$, they follow from Proposition \ref{prop:ovhol}, Theorem \ref{thm:ovhol-dev} and compatibility with Frobenius pullback.
	\end{remark}
	
	\begin{remark}
		When $X$ is a realizable $k$-scheme, i.e., there is a frame $(X, Y, \Pf)$ with $Y$ proper $k$ and $\Pf$ projective smooth over $\Pf$, the triangulated subcategory $D^b_{\oh, F}$ defined here is the same as $D^b_{\textrm{hol}, F}(X/K)$ defined in \cite[\S 4.2]{lazda2023rigidification}. This is because that we have a closed immersion $(X, Y, \Pf) \hookrightarrow (X, Y', \Pf') = (X, Y \times_{\Spec k} k[[t]], \Pf \times_{\Spf \V} \Spf \V[[t]])$ and this induces the equivalence $\D^b_{\oh}(X, Y, \Pf) \simto D^b_{\oh}(X, Y', \Pf')$ similarly as in Lemma \ref{lem:indep-frame}. Also six functors coincide.
	\end{remark}
	
	\begin{theorem} \label{thm:sp-compatibility}
		Let $X$ be a realizable scheme over $k[[t]]$ and $f : X' \to X$ be a morphism of realizable schemes over $k[[t]]$. For $E, F \in \Isoc_F(X/K)$, we have
		$$ \sp_{X, +}(E^\vee) \simeq \Db_X \sp_{X, +}(E),$$
		$$ \sp_{X, +}(E \otimes F) \simeq \sp_{X, +}(E) \, \wt{\otimes} \, \sp_{X, +}(F)[\delta_X]$$
		and
		$$ \sp_{X', +}f^*_K(E) \simeq f^! \sp_{X, +}(E)[\delta_{X'} - \delta_X].$$ 
	\end{theorem}
	\begin{proof}
		By Proposition \ref{prop:sp-dual}, \ref{prop:sp-pullback} and \ref{prop:sp-tensor}.
	\end{proof}
	
	\begin{remark}
		The shifts of the compatibility in Theorem \ref{thm:sp-compatibility} should disappear by introducing the dual constructible $t$-structure on $D^b_{\oh, F}$ and a $t$-exact functor $\sp^! : D^b_{\oh, F}(X/K) \to D^b_{\mathrm{cons}, F}(X/K)$ as in \cite{lazda2023rigidification}. This will be studied in the future.
	\end{remark}
	
	\bibliographystyle{alpha}
	\bibliography{main}

@article{caro2012overholonomicity,
	title={Overholonomicity of overconvergent {$F$}-isocrystals over smooth varieties},
	author={Caro, Daniel and Tsuzuki, Nobuo},
	journal={Annals of Mathematics},
	pages={747--813},
	year={2012},
	publisher={JSTOR}
}

@article{Caro2011,
	author    = {Caro, Daniel},
	title     = {Pleine fid{\'e}lit{\'e} sans structure de Frobenius et isocristaux partiellement surconvergents},
	journal   = {Mathematische Annalen},
	volume    = {349},
	number    = {4},
	pages     = {747--805},
	year      = {2011},
	doi       = {10.1007/s00208-010-0539-x},
	url       = {https://link.springer.com/article/10.1007/s00208-010-0539-x}
}

@article{liu2026semistable,
	title={Semistable Reduction Theorem for Overconvergent {$F$}-isocrystals over Laurent Series Fields},
	author={Liu, Yuanmin},
	journal={arXiv preprint arXiv:2604.17799},
	year={2026}
}

@article{caro2022logarithmic,
	author = {Caro, Daniel and Vauclair, David},
	title = {Logarithmic $p$-bases and arithmetical differential modules},
	journal = {Annales de l'Institut Fourier},
	pages = {1011--1096},
	year = {2022},
	publisher = {Association des Annales de l{\textquoteright}institut Fourier},
	volume = {72},
	number = {3},
	doi = {10.5802/aif.3508},
	language = {en},
	url = {https://aif.centre-mersenne.org/articles/10.5802/aif.3508/}
}

@article{caro2021arithmetic,
	title={Arithmetic {$\mathcal{D}$}-modules over Laurent series fields: absolute case},
	author={Caro, Daniel},
	journal={arXiv preprint arXiv:2103.10362},
	year={2021}
}

@article{abe2018theory,
	title={Theory of weights in $p$-adic cohomology},
	author={Abe, Tomoyuki and Caro, Daniel},
	journal={American Journal of Mathematics},
	volume={140},
	number={4},
	pages={879--975},
	year={2018},
	publisher={Johns Hopkins University Press}
}

@article{abe2025proper,
	title={Proper pushforwards on analytic adic spaces},
	author={Abe, Tomoyuki and Lazda, Christopher},
	journal={The Quarterly Journal of Mathematics},
	volume={76},
	number={1},
	pages={147--183},
	year={2025},
	publisher={Oxford University Press UK}
}

@article{caro2006,
	title={D\'{e}vissages des {$F$}-complexes de $\mathcal{D}$-modules arithm\'{e}tiques en {$F$}-isocristaux surconvergents},
	author={Caro, Daniel},
	journal={Inventiones Mathematicae},
	volume={166},
	pages={397–456},
	year={2006},
	url={https://doi.org/10.1007/s00222-006-0517-9}
}

@article{berthelotDI,
	author = {Berthelot, Pierre},
	title = {{$\mathcal{D}$}-modules arithm\'etiques. {I.} {Op\'erateurs} diff\'erentiels de niveau fini},
	journal = {Annales scientifiques de l'\'Ecole Normale Sup\'erieure},
	pages = {185--272},
	year = {1996},
	publisher = {Elsevier},
	volume = {4e s{\'e}rie, 29},
	number = {2},
	doi = {10.24033/asens.1739},
	mrnumber = {97b:14019},
	zbl = {0886.14004},
	language = {fr},
	url = {https://www.numdam.org/articles/10.24033/asens.1739/}
}

@article{berthelotDII,
	author = {Berthelot, Pierre},
	title = {{$\mathcal{D}$}-modules arithm\'etiques. {II} : descente par {F}robenius},
	series = {M\'emoires de la Soci\'et\'e Math\'ematique de France},
	year = {2000},
	publisher = {Soci\'et\'e math\'ematique de France},
	number = {81},
	doi = {10.24033/msmf.394},
	mrnumber = {2001k:14043},
	zbl = {0948.14017},
	language = {fr},
	url = {https://www.numdam.org/item/MSMF_2000_2_81__1_0/}
}

@article{lazda2016rigid,
	title={Rigid cohomology over {L}aurent series fields},
	author={Lazda, Christopher and P{\'a}l, Ambrus},
	volume={21},
	year={2016},
	publisher={Springer}
}

@article{caro2019arithmetic,
	title={Arithmetic $\mathcal{D}$-modules over algebraic varieties of characteristic {$p > 0$}},
	author={Caro, Daniel},
	journal={arXiv preprint arXiv:1909.08432},
	year={2019}
}

@misc{caro2025arith,
	title={Arithmetic $\mathcal{D}$-modules},
	author={Caro, Daniel},
	url={https://carod.users.lmno.cnrs.fr/Arithmetic%20D-modules.pdf},
	note={available at \url{https://carod.users.lmno.cnrs.fr/Arithmetic%20D-modules.pdf}, accessed at Oct 1, 2026}
}

@article{lazda2023rigidification,
	title={Rigidification of arithmetic $\mathscr{D}$-modules and an overconvergent Riemann-Hilbert correspondence},
	author={Lazda, Christopher},
	journal={arXiv preprint arXiv:2304.07181},
	year={2023}
}

@unpublished{berthelot1996cohomologie,
	author    = {Pierre Berthelot},
	title     = {Cohomologie rigide et cohomologie rigide à support propre, premi{\`e}re partie},
	note      = {Pr{\'e}publication IRMAR 96-03},
	school    = {Universit{\'e} de Rennes 1},
	year      = {1996},
	url       = {https://www.wstein.org/people/berthelo/publis/Cohomologie_Rigide_I.pdf}
}

@article{de1996smoothness,
	title={Smoothness, semi-stability and alterations},
	author={de Jong, Aise Johan},
	journal={Publications Math{\'e}matiques de l'IH{\'E}S},
	volume={83},
	pages={51--93},
	year={1996}
}

@article{shiho2002crystallineII,
	title={Crystalline fundamental groups. II. Log convergent cohomology and rigid cohomology},
	author={Shiho, Atsushi},
	journal={J. Math. Sci. Univ. Tokyo},
	volume={9},
	number={1},
	pages={1--163},
	year={2002}
}

@article{andre2002HasseArf,
	author    = {Andr{\'e}, Yves},
	title     = {Filtrations de type {Hasse--Arf} et monodromie $p$-adique},
	journal   = {Inventiones Mathematicae},
	volume    = {148},
	number    = {2},
	pages     = {285--317},
	year      = {2002},
	doi       = {10.1007/s002220100207},
}

@article{mebkhout2002analogue,
	title={Analogue $p$-adique du th{\'e}oreme de Turrittin et le th{\'e}oreme de la monodromie $p$-adique},
	author={Mebkhout, Zoghman},
	journal={Inventiones mathematicae},
	volume={148},
	number={2},
	pages={319--351},
	year={2002},
	publisher={Springer}
}

@article{kedlaya2004p,
	title={A $p$-adic local monodromy theorem},
	author={Kedlaya, Kiran S},
	journal={Annals of mathematics},
	pages={93--184},
	year={2004},
	publisher={JSTOR}
}

@article{kedlaya2007semistable,
	title={Semistable reduction for overconvergent {F}-isocrystals {I}: Unipotence and logarithmic extensions},
	author={Kedlaya, Kiran},
	journal={Compositio Mathematica},
	volume={143},
	number={5},
	pages={1164--1212},
	year={2007},
	publisher={London Mathematical Society}
}

@article{kedlaya2008semistable,
	title={Semistable reduction for overconvergent {F}-isocrystals, {II}: A valuation-theoretic approach},
	author={Kedlaya, Kiran},
	journal={Compositio Mathematica},
	volume={144},
	number={3},
	pages={657--672},
	year={2008},
	publisher={London Mathematical Society}
}

@article{kedlaya2009semistable,
	title={Semistable reduction for overconvergent {F}-isocrystals, {III}: Local semistable reduction at monomial valuations}, 
	volume={145},
	DOI={10.1112/S0010437X08003783},
	number={1},
	journal={Compositio Mathematica},
	author={Kedlaya, Kiran},
	year={2009},
	pages={143--172}
}

@article{kedlaya2011semistable,
	title={Semistable reduction for overconvergent {F}-isocrystals, {IV}: Local semistable reduction at nonmonomial valuations},
	author={Kedlaya, Kiran},
	journal={Compositio Mathematica},
	volume={147},
	number={2},
	pages={467--523},
	year={2011},
	publisher={London Mathematical Society}
}

@article{kato1991explicit,
	title={The explicit reciprocity law and the cohomology of {F}ontaine-{M}essing},
	author={Kato, Kazuya},
	journal={Bulletin de la Soci{\'e}t{\'e} Math{\'e}matique de France},
	volume={119},
	number={4},
	pages={397--441},
	year={1991}
}

@article{shiho2007relative,
	title={Relative log convergent cohomology and relative rigid cohomology {II}},
	author={Shiho, Atsushi},
	journal={arXiv preprint arXiv:0707.1743},
	year={2007}
}
\end{document}